\documentclass[reqno,10pt,a4paper]{amsart}
\usepackage{graphicx,color,amssymb,latexsym}
\usepackage{mathrsfs}

\usepackage{lscape}
\newtheorem{corollary}{\bf Corollary}
\newtheorem{remark}{Remark}

\newtheorem{theorem}{\bf Theorem}[section]

\numberwithin{equation}{section} \theoremstyle{plain}
\theoremstyle{definition}

\newtheorem{example}{Example}[section]

\DeclareMathOperator*{\argmax}{argmax}
\DeclareMathOperator*{\cov}{cov} \DeclareMathOperator{\diag}{diag}
\DeclareMathOperator*{\trace}{trace}
\DeclareMathOperator*{\sign}{sign}
\newcommand{\lef}{\langle\hskip -1.8pt \langle}
\newcommand{\rig}{\rangle\hskip -1.8pt \rangle}

\input cyracc.def

\begin{document}

\title[Examples of moderate deviation principle ]
{Examples of moderate deviation principle for diffusion processes}

\author{A. Guillin}
\address{CEREMADE, Universit\'e Paris Dauphine and TSI, Ecole nationale des
Telecommunications}
\email{guillin@ceremade.dauphine.fr}

\author{R. Liptser}
\address{Electrical Engineering Systems,
Tel Aviv University, 69978 - Ramat Aviv, Tel Aviv, Israel}
\email{liptser@eng.tau.ac.il} 


\date{December 25, 2004.}

\maketitle
\begin{abstract}
Taking into account some likeness of moderate deviations (MD) and
central limit theorems (CLT), we develop an approach, which made a
good showing in CLT, for MD analysis of a family
$$ S^\kappa_t=\frac{1}{t^\kappa}\int_0^tH(X_s)ds, \
t\to\infty
$$
for an ergodic diffusion process $X_t$ under
$0.5<\kappa<1$ and appropriate $H$. We mean a decomposition with
``corrector'':
$$
\frac{1}{t^\kappa}\int_0^tH(X_s)ds={\rm
corrector}+\frac{1}{t^\kappa}\underbrace{M_t}_{\rm martingale}.
$$
and show that, as in the CLT analysis, the corrector is negligible but
in the MD scale, and the main contribution in the MD brings the
family ``$ \frac{1}{t^\kappa}M_t, \ t\to\infty. $'' Starting from
Bayer and Freidlin, \cite{BF}, and finishing by Wu's papers
\cite{Wu1}-\cite{WuH}, in the MD study Laplace's transform
dominates. In the paper, we replace the Laplace technique by one,
admitting to give the conditions, providing the MD, in terms of
``drift-diffusion'' parameters and $H$. However, a verification of
these conditions heavily depends on a specificity of a diffusion
model. That is why the paper is named ``Examples ...''.
\end{abstract}

\section{Introduction}
\label{sec-1}

In this paper, we study the moderate deviation principle (in
short: MDP) for a family $(S^\kappa_t)_{t\to\infty}$, $\kappa\in
\big(\frac{1}{2},1\big)$:
\begin{equation*}
S^\kappa_t=\frac{1}{t^\kappa}\int_0^tH(X_s)ds,
\end{equation*}
where $X=(X_t)_{t\ge 0}$ is an ergodic diffusion process ($X_t\in
\mathbb{R}^d$, $d\ge 1$) (with the unique invariant measure $\mu(dz)$,
obeying the density $p(z)$ relative to Lebesgue measure over
$\mathbb{R}^d$.

The function $H:\mathbb{R}^d\to\mathbb{R}^q$ is assumed to be
integrable relative to $\mu(dz)$ and has zero barycenter
\begin{equation}\label{nul}
\int_{\mathbb{R}^d}H(z)p(z)dz=0.
\end{equation}
We restrict ourselves by consideration of  the strong (unique)
solution of It\^o's equation
\begin{equation}\label{ito00}
dX_t=b(X_t)dt+\sigma(X_t)dW_t
\end{equation}
generated by a standard vector-valued Wiener process
$W=(W_t)_{t\ge 0}$ and subject to a fixed initial point, $X_0=x$.
We also include into the consideration a linear version of
\eqref{ito00} (here $A,B$ are matrices):
\begin{equation}\label{itolin}
dX_t=AX_tdt+BdW_t
\end{equation}
being popular in engineering.

In a nonlinear case, we use Veretennikov - Khasminskii's condition
(see, \cite{Kh} and \cite{Ver98}): for some positive numbers $r$, $C$ and
$\alpha$, (here $\lef\cdot\rig$ denotes the inner product)
$$
\lef z, b(z)\rig\le -r\|z\|^{1+\alpha}, \ \|z\|>C
$$
and assume that the diffusion matrix $ a(x)=\sigma\sigma^*(x) $ is
nonsingular and bounded. In a linear case, proper assumptions are
given in terms of the pair $(A,B)$:

1) eigenvalues of $A$ have negative real parts;

2) $(A,B)$ satisfies Kalman's controllability condition from
\cite{O77},

\hskip .2in i.e. a singularity of $a(x)\equiv BB^*$ is
permissible.

\noindent
For the  MDP analysis, we apply well known method employed  for the
central limit theorem (in short CLT) proof of a family
$$
\Big(\frac{1}{\sqrt{t}}\int_0^tH(X_s)ds\Big)_{t\to\infty}
$$
(see, e.g. Papanicolaou, Stroock and Varadhan \cite{PSV}, Ethier
and Kurtz \cite{EK}, Bhattacharya \cite{Br}, Pardoux and
Veretennikov \cite{PVI}, \cite{PVII} and citations therein, see
also Ch. 9, \S 3 in \cite{LSMar}) based on a decomposition with
corrector:
$$
\frac{1}{\sqrt{t}}\int_0^tH(X_s)ds=\frac{1}{\sqrt{t}}\underbrace{[U(x)-U(X_t)]}_
{\rm corrector}+\frac{1}{\sqrt{t}}\underbrace{M_t}_{\rm
martingale},
$$
where
\begin{equation*}\label{1.5UU}
U(x)=\int_0^\infty\int_{\mathbb{R}^d}H(y)P^{(t)}_x(dy)dt,
\end{equation*}
$P^{(t)}_x$ is the transition probability kernel of $X$, and $M_t$
is a continuous martingale with the variation process $\langle M\rangle_t$. In the above mentioned
papers, the corrector is negligible in a sense
$$
\frac{1}{\sqrt{t}}[U(x)-U(X_t)]\xrightarrow[t\to\infty]{\rm
prob.}0
$$
and the main contribution to a limit distribution brings
$\frac{1}{\sqrt{t}}M_t$. It is well known (see, e.g. Ch. 5 in
\cite{LSMar}) the following implication: with nonnegative definite
matrix
$$
\frac{1}{t}\langle M\rangle_t\xrightarrow[t\to\infty]{\rm prob.}Q
\Rightarrow Ee^{\lef\lambda,\frac{1}{\sqrt{t}}M_t\rig}
\xrightarrow[t\to\infty]{} e^{-\frac{1}{2}\lef \lambda,
Q\lambda\rig},\quad \forall \ \lambda\in\mathbb{R}^q,
$$
where
$$
Q=\int_0^\infty\int_{\mathbb{R}^d}
\big[(P^{(t)}_zH)H^*(z)+(P^{(t)}_zH)^*H(z)\big]p(z)dzdt,
$$

Summarizing these remarks, we may claim that the CLT holds provided
that $U(x)$ and $Q$ exist and for any $\varepsilon>0$
\begin{eqnarray*}
&&\lim\limits_{t\to\infty}P\big(|U(x)-U(X_t)|>\sqrt{t}\varepsilon\big)=0
\\
&&\lim\limits_{t\to\infty}P\big(|\langle
M\rangle_t-tQ|>t\varepsilon\big)=0.
\end{eqnarray*}

We develop the same method for MDP analysis.
Replacing $\frac{1}{\sqrt{t}}$ by $\frac{1}{t^\kappa}$, we keep
the CLT framework with the same $U(x)$, $M_t$ and $Q$, i.e.,
\begin{equation*}
\frac{1}{t^\kappa}\int_0^tH(X_s)ds=\frac{1}{t^\kappa}\underbrace{[U(x)-U(X_t)]}_
{\rm corrector}+\frac{1}{t^\kappa}\underbrace{M_t}_{\rm
martingale}
\end{equation*}
and claim that (Theorem \ref{theo-main}) the MDP holds, with the
rate of speed
$$
\varrho(t)= \frac{1}{t^{2\kappa-1}}
$$
provided that $U(x)$ and $Q$ exist and for any $\varepsilon>0$
\begin{equation}\label{1.6ux}
\begin{aligned}
& \lim\limits_{t\to\infty}\varrho(t)\log
P\big(|U(x)-U(X_t)|>t^\kappa\varepsilon\big)= -\infty
\\
&\lim\limits_{t\to\infty}\varrho(t)\log P\big(|\langle
M\rangle_t-tQ|>t\varepsilon\big)= -\infty.
\end{aligned}
\end{equation}
A choice of $\varrho(t)$ is imposed by $\frac{1}{t^\kappa}$. As in
the CLT proof, the corrector negligibility is required but
exponentially fast with the rate of speed $\varrho(t)$. The
main contribution in the MDP brings the family
$\big(\frac{1}{t^\kappa}M_t\big)_{t\to\infty}$.

Most probably, Dembo, \cite{D}, was one of the first who introduced  a condition
of \eqref{1.6ux} (second) type. We  found in Puhalskii, \cite{puh1} (Theorem 2.3)
and \cite{P3}, \cite{puh2} that, in our setting with nonsingular (!) matrix $Q$,
\eqref{1.6ux} provides MDP for the family
$\big(\frac{1}{t^\kappa}M_t\big)_{t\to\infty}$ with the rate of
speed $\varrho(t)$ and the rate function
$$
J(Y)=\frac{1}{2}\|Y\|^2_{Q^{-1}}, \ Y\in\mathbb{R}^q.
$$
We prove
in Theorem \ref{theo-main} that the same statement remains valid
for a singular $Q$ too with the rate function

\begin{equation*}
J(Y)=\begin{cases} \frac{1}{2}\|Y\|^2_{Q^\oplus},& \ Y=QQ^\oplus
Y,
\\
\infty, &\text{otherwise},
\end{cases}
\end{equation*}
where $Q^\oplus$ is the Moore-Penrose pseudoinverse matrix (see,
Albert, \cite{Albert}).

It would be noted that seeming simplicity of \eqref{1.6ux} is
delusive with the exception of the eigenvalue gap case (in short EG, see
Gong and Wu, \cite{GongWu5}) for $P^{(t)}_x$ (a corresponding scenario
can be found in
\cite{DJL}). Unfortunately, the EG fails for diffusion processes.
For instance, under $P^{(t)}_x$ associated with Ornstein-Uhlenbeck's process
$$
dX_t=-X_tdt+dW_t
$$
having $\big(0,\frac{1}{2}\big)$-Gaussian
invariant measure $\mu$, if EG were valid, then for bounded centered $H$
$$
|P^{(t)}_xH|\le \text{cons.}e^{-\lambda t}, \
\forall \ t\ge 0, \ \exists \ \lambda>0.
$$
However, direct computations show that for $H(x)=\sign(x)$ and sufficiently large
$|x|$, we have
$
|P^{(t)}_xH|dt\le \upsilon(x)e^{-\lambda t}
$
where $\upsilon(x)$ is a positive function,
$\upsilon(x)<\infty$ over $\mathbb{R}^d$ and
$\upsilon(x)\to\infty$ with $|x|\to\infty$. The condition of this
type: for any bounded and measurable $H$
$$
|P^{(t)}_xH-\mu H|\le\upsilon(x)e^{-\lambda t}
$$
describes the geometric ergodicity (see, Down, Meyn and
Tweedie, \cite{DMT} and citations therein). The geometric ergodicity is
a helpful tool for the verification of $U(x)$ and $Q$ existence and even
for the first part of \eqref{1.6ux} verification, although, a crude choice
of $\upsilon(x)$, say $\upsilon(x)\asymp |x|^m, m>2$, may to render this
verification impossible (CLT analysis is not so sensitive to a choice of $\upsilon$).
The second part of \eqref{1.6ux} verification is very sensitive to properties of $U$,
owing to $\langle M\rangle_t=\int_0^t\nabla^* U(X_s)(a(X_s)\nabla U(X_s)ds$,
so that, the geometric ergodicity framework is not a ``foreground'' tool.
Following Pardoux and Veretennikov, \cite{PVI},  we combine  a property of $H$
with a polynomial ergodicity $|P^{(t)}_xH-\mu H|\le \frac{\upsilon(x)}
{(1+t)^\gamma}$, $\gamma>1$ with $H$-depending $\upsilon$ admitting
an effective verification of \eqref{1.6ux}.
In this connection, we mention here some result (see, Theorem \ref{theo-a.1}),
in Appendix, interesting by itself, which is helpful in \eqref{1.6ux} verification.
Let $X$ be a diffusion
process with the generator $\mathscr{L}$ and $V(x)$ is Lyapunov's
function belonging to the range of definition of $\mathscr{L}$.
Then,
$
N_t=V(X_t)-V(x_0)-\int_0^t\mathscr{L}V(X_s)ds
$
is a continuous martingale and denote by $\langle N\rangle_t$ its
variation process. Assume:
$$
\mathscr{L}V\le -cV^\ell+\mathfrak{c}, \ \exists \ q>0 \ \text{and} \
\langle N\rangle_t\le \int_0^t\mathbf{c}\big(1+V^r(X_s)\big)ds, \
\exists \ r\le \ell.
$$
Then, for any $\varepsilon>0$ and sufficiently large number $n$
$$
\begin{aligned}
&
 \lim\limits_{t\to\infty}\varrho(t)\log
P\big(V(X_t)>t^{2\kappa}\varepsilon\big) =-\infty,
\\
& \lim\limits_{t\to\infty}\varrho(t)\log P\Big(\int_0^tV^\ell(X_s)ds
>tn\Big)
=-\infty.
\end{aligned}
$$

Our method of the MDP analysis differs from Wu \cite{Wu1} -
\cite{WuH} where the Laplace transform technique dominates, or
Guillin \cite{G1}, \cite{G2} based on discrete time approximation
and Markov chains. In our approach, we deal with the
above-mentioned Puhalskii's results obtained with the help of, so
called, stochastic exponential as an alternative to Laplace's transform
technique (see, e.g. \cite{DJL} for more detailed explanation in
the discrete time case).

The paper is organized as follows. In Section \ref{sec-00}, all
notations are given and Theorem \ref{theo-main}, generalized
Puhalskii's for singular $Q$, is formulated and proved. In Section
\ref{sec-3new}, all results and examples are presented focusing on the existence
and properties of the corrector and martingale variation process. The proofs are
gathered in Section \ref{sec-4}. A simple example showing how the MDP may help
in a statistical inference (for more information on statistical applications
see, Inglot and Kallenberg, \cite{IngKaL}) is given in Section
\ref{sec-5}.
The technical tools are gathered in Appendix \ref{sec-A}.

\section{Preliminaries}
\label{sec-00}

We fix the following notations and assumptions which are in force through the
paper.
The random process $X=(X_t)_{t\ge 0}$ is defined on some
stochastic basis $(\Omega,\mathscr{F},\mathbf{F}=(\mathscr{F}_t)
_{t\ge },P)$ satisfying the usual conditions.

\bigskip

-  $\|\cdot\|$, $|\cdot|$, and $\lef\cdot,\cdot\rig$ are
Euclidean's and $\mathbb{L}$ norms respectively in $\mathbb{R}^d$ and the

\hskip .1in inner product.

\smallskip
- $^*$ is transposition symbol.

\smallskip
- $ a(z):=\sigma\sigma^*(z). $

\smallskip
- $c, \mathfrak{c}, \mathbf{c}\in\mathbb{R}_+$, $\ldots,$ are generic
constants.

\smallskip
- $P^{(t)}_x(dy)$ is the transition probability kernel of $X$.

\smallskip
- $E_x$ denotes the expectation relative to  $P^{(t)}_x(dy)$.

\smallskip
- $\mu(dz)$ is the invariant measure.

\smallskip
- $ \mathscr{L}=
\frac{1}{2}\sum_{i,j=1}^{d}a_{ij}(z)\frac{\partial^2} {\partial
z_i\partial z_j} +\sum_{i=1}^{d}b_i(z)\frac{\partial}{\partial
z_i} $ is the generator of $X$.

\smallskip
- $(\mathscr{F}^X_t)_{t\ge 0}$ is the filtration, with the general
conditions, generated by $(X_t)$.

\smallskip
- $\langle L\rangle_t$ - is the variation process of a continuous
martingale $(L_t)_{t\ge 0}.$

\smallskip
- $\nabla f(x)$ is the gradient of $f(x)$ (row vector).

\smallskip
- $\rho$ is Euclidean's metric in $\mathbb{R}^d$.

\smallskip
- $\varrho(t)=\frac{1}{t^{2\kappa-1}}$.

\smallskip
- $\mathrm{I}$ denotes the identical matrix of an appropriate
size.

\smallskip
- ``$>$'', ``$\ge $'' denote also the standard inequalities for
nonnegative definite matri-

\hskip .1in ces.

\bigskip
As was mentioned in Introduction, the existence of
\begin{eqnarray}\label{QQQQQ}
Q&=&\int_0^\infty\int_{\mathbb{R}^d}
\Big[(P^{(t)}_zH)H^*(z)+(P^{(t)}_zH)^*H(z)\Big]dtp(z)dz,
\\
U(x)&=&\int_0^\infty\int_{\mathbb{R}^d}H(y)P^{(t)}_x(dy)dt
\label{1.5U}
\end{eqnarray}
is required. We emphasize that
\begin{equation*}
M_t=U(X_t)-U(x)+\int_0^tH(X_s)ds
\end{equation*}
is the martingale relative to $(\mathscr{F}^X_t)_{t\ge 0}$.

\smallskip
The theorem below is a ``master-key''  for MDP analysis.

\begin{theorem}\label{theo-main}
For any $x\in \mathbb{R}^d$ and any $\varepsilon>0$, assume

{\rm (i)} $ \lim\limits_{t\to\infty}\varrho(t)\log
P\big(|U(x)-U(X_t)|>t^\kappa\varepsilon\big)= -\infty $

\smallskip
{\rm (ii)} $ \lim\limits_{t\to\infty}\varrho(t)\log P\big(|\langle
M\rangle_t-tQ|>t\varepsilon\big)= -\infty $

\smallskip
\noindent Then, the family $(S^\kappa_t)_{t\to\infty}$ obeys the
MDP in $(\mathbb{R}^q,\rho)$ with the rate of speed $\varrho(t)$
and the rate function
\begin{equation}\label{1.7e}
J(Y)=\begin{cases} \frac{1}{2}\|Y\|^2_{Q^\oplus},& \ Y=QQ^\oplus
Y,
\\
\infty, &\text{otherwise},
\end{cases}
\end{equation}
where $Q^\oplus$ is the Moore-Penrose pseudoinverse matrix {\rm
(}see, Albert, \cite{Albert}{\rm )}.
\end{theorem}
\begin{proof}
From the definition of $M_t$, it follows that
$$
S^\kappa_t=\frac{1}{t^\kappa}[U(x)-U(X_t)]+\frac{1}{t^\kappa}M_t.
$$

(i) provides the negligibility of
$\big(\frac{1}{t^\kappa}[U(x)-U(X_t)]\big)_{t\to\infty}$ in $\varrho$-MDP scale.

(ii) provides $\varrho$-MDP, , under positive definite matrix, with
the rate function
\begin{equation*}
J(Y)=\frac{1}{2}\|Y\|^2_{Q^{-1}}
\end{equation*}
for the family $\big(\frac{1}{t^\kappa}M_t\big)_{t\to\infty}$
(due to result similar to Puhalskii, \cite{puh1} (Theorem 2.3) and \cite{puh2}).

\medskip
If $Q$ is nonnegative definite only, the above result is no longer valid.
This remark necessitates to turn to the general approach in large deviation
analysis adapted to our setting. The family
$\big(\frac{1}{t^\kappa}M_t\big)_{t\to\infty}$ is said to obey the
large deviation principle (in our terminology: MDP) with the rate
of speed $\varrho(t)$ and some (good) rate function $J(Y),Y\in
\mathbb{R}^q$, provided that this family is
$\varrho$-exponentially tight in $(\mathbb{R}^q,\rho)$:
\begin{equation}\label{expt}
\lim_{K\to \infty}\varlimsup_{t\to \infty}\varrho(t)\log P
\Big(\Big|\frac{1}{t^\kappa}M_t\Big|>K\Big)=-\infty
\end{equation}
and obeys $(\varrho,J)$-local large deviation principle with the
rate function $J(Y)$: for any $Y\in \mathbb{R}^q$
\begin{equation}\label{2.21sing1}
\begin{aligned}
& \varlimsup_{\delta\to 0}\varlimsup_{t\to \infty}\varrho(t)\log P
\Big(\Big|\frac{1}{t^\kappa}M_t-Y\Big|\le\delta\Big)\le -J(Y)
\\
& \varliminf_{\delta\to 0}\varliminf_{t\to \infty}\varrho(t)\log P
\Big(\Big|\frac{1}{t^\kappa}M_t-Y\Big|\le \delta\Big)\ge -J(Y).
\end{aligned}
\end{equation}
A direct verification of \eqref{expt} and \eqref{2.21sing1} would
be difficult. So, it is reasonable to verify \eqref{expt}
by applying the following regularization procedure. We introduce a new family
$\big(\frac{1}{t^\kappa}M^\gamma_t\big)_{t\to\infty}$ with
$$
M^\gamma_t=M_t+\sqrt{\gamma}dW'_t,
$$
where $\gamma$ is a positive number and $W'_t (\in \mathbb{R}^q)$
is a standard Wiener process independent of $M_t$. The random
process $M^\gamma_t$ is continuous martingale with
$$
\langle M^\gamma\rangle_t=\langle M\rangle_t+\gamma
\mathrm{I}t,
$$
where $\mathrm{I}=\mathrm{I}_{q\times q}$.
For the family $\big(\frac{1}{t^\kappa}M^\gamma_t\big)_{t\to\infty}$,
(ii) reads as:
\[
\lim_{t\to\infty}\varrho(t)\log P\Big(\Big| \frac{1}{t}\langle
M^\delta\rangle_t-Q_\gamma\Big|>\varepsilon\Big)=-\infty,
\]
where $Q_\gamma=Q+\gamma\mathrm{I}$.
Since $Q_\gamma$ is the nonsingular matrix, the family
$\big(\frac{1}{t^\kappa}M^\gamma_t\big)_{t\to\infty}$ obeys $(\varrho,J_\gamma)$-MDP,
where
$
J_\gamma(Y)=\frac{1}{2}\|Y\|^2_{Q^{-1}_\gamma}
$.

Now, we apply the basic  Puhalskii theorem from \cite{P1} which,
being adapted to our case, states that the family
$\big(\frac{1}{t^\kappa}M^\gamma_t\big)_{t\to\infty}$ is
$\varrho(t)$-exponentially tight, in $(\mathbb{R}^q,\rho)$:
\begin{equation}\label{2.2q}
\varlimsup_{K\to \infty}\varlimsup_{t\to \infty}\varrho(t)\log P
\Big(\Big|\frac{1}{t^\kappa}M^\gamma_t\Big|>K\Big) =-\infty,
\end{equation}
and obeys $(\varrho,J_\gamma)$-local deviation principle:
\begin{equation}\label{2.21song1}
\begin{aligned}
& \varlimsup_{\delta\to 0}\varlimsup_{t\to \infty}\varrho(t)\log P
\Big(\Big|\frac{1}{t^\kappa}M^\gamma_t-Y\Big|\le\delta\Big) \le
-J_\gamma (Y)
\\
& \varliminf_{\delta\to 0}\varliminf_{t\to \infty}\varrho(t)\log P
\Big(\Big|\frac{1}{t^\kappa}M^\gamma_t-Y\Big|\le\delta\Big) \ge
-J_\gamma(Y).
\end{aligned}
\end{equation}
Obviously, \eqref{2.2q} and \eqref{2.21song1} imply \eqref{expt} and
\eqref{2.21sing1} provided that
\begin{equation}\label{2.24end1}
\lim_{\delta \to 0}\varlimsup_{t\to \infty}\varrho(t)P
\Big(\Big|\frac{\sqrt{\gamma}}{t^\kappa}W'_t\Big|\ge
\eta\Big)=-\infty, \quad \ \forall \ \eta>0
\end{equation}
and
\begin{equation}\label{2.23+1}
\lim_{\gamma \to 0}J_\gamma (V)=
\begin{cases}
\frac{1}{2}\|Y\|^2_{Q^\oplus}, & Q^\oplus QY=Y
\\
\infty, & \text{otherwise}.
\end{cases}
\end{equation}

\eqref{2.24end1} holds true, since the family
$\frac{\sqrt{\gamma}}{t^\kappa}W'_t$ obeys the $\varrho$-MDP with the rate
function $
\frac{1}{2\gamma}\|Y\|^2$, so that,
$$
\varlimsup_{t\to \infty}\frac{1}{t^{2\kappa-1}}{}P
\Big(\Big\|\frac{\sqrt{\gamma}}{t^\kappa}W'_t\Big\|\ge
\eta\Big)\le -\inf_{\{Y:\|Y\|\ge \frac{\eta}{\sqrt{\gamma}}\}}
\frac{1}{2}\|Y\|^2=-\frac{\eta^2}{2\gamma}\xrightarrow[\gamma\to
0]{} -\infty.
$$

\eqref{2.23+1} is verified with an utilization of the pseudoinverse matrix properties.
Let $T$ be an orthogonal matrix transforming $Q$ to
the diagonal form: $ \diag(Q)=T^*QT.$ Due to
$$
2J_\gamma (Y)=Y^*\big[\gamma\mathrm{I}+Q\big]^{-1}Y=
Y^*T\big[\gamma\mathrm{I}+\diag(Q)\big]^{-1} T^*Y,
$$
for $Y=Q^\oplus QY$ we have (recall that $Q^\oplus
QQ^\oplus=Q^\oplus$, see \cite{Albert})
$$
\begin{aligned}
2J_\gamma (Y)&=Y^*Q^\oplus QT\big[\gamma
\mathrm{I}+\diag(Q)\big]^{-1}T^*Y
\\
&=Y^*Q^\oplus TT^* QT\big[\gamma\mathrm{I}+\diag(Q)\big]^{-1}
T^*Y
\\
&=Y^*Q^\oplus T\diag(Q)\big[\gamma \mathrm{I}+\diag(Q))^{-1}T^*Y
\\
&\xrightarrow[\gamma \to 0]{}Y^*Q^\oplus
T\diag(Q)\diag^\oplus(Q)T^*Y
\\
&=Y^*Q^\oplus T\diag(Q)T^*T\diag^\oplus(Q)T^*Y
\\
&=Y^*Q^\oplus QQ^\oplus Y=Y^*Q^\oplus Y=\|Y\|^2_{Q^\oplus}=2J(Y).
\end{aligned}
$$

For $Y\ne Q^\oplus QY$, $\lim_{\gamma\to 0}J_\gamma(Y)=\infty$.
\end{proof}

\section{Main results}\label{sec-3new}
\vskip .2in
\subsection{Nonlinear model, I}\label{sec-01}

$X_t$ solves \eqref{ito00} subject to $X_0=x$.

\bigskip
$\mathbf{(A_b)}$ \ \ $b$ is locally Lipschitz continuous; for some
$\alpha\ge 1$ and $C>0$ there

\hskip .48in exists $\mathfrak{r}>0$, depending on $\alpha,C$,
such that
$$
\lef z, b(z)\rig\le -\mathfrak{r}\|z\|^{1+\alpha}, \ \|z\|>C.
$$

\medskip
$\mathbf{(A_{\sigma, a})}$ $\sigma$ is Lipschitz continuous; for
some $\Lambda>\lambda>0$
$$
\lambda \mathrm{I}\leq a(z)\leq \Lambda \mathrm{I}.
$$

From Pardoux and Veretennikov \cite{PVI}, it follows that, under
$\mathbf{(A_b)}$ and $\mathbf{(A_{\sigma, a})}$, the diffusion process $X$
is ergodic with the unique invariant measure $\mu(dz)$ possessing
a density $p(z)$ relative to $dz$. Moreover, for $\alpha>1$ and
any $\beta<0$
$$
\int_{\mathbb{R}^d}(1+\|z\|)^{\alpha-1+\beta}p(z)dz<\infty.
$$

\medskip
$\mathbf{(A_H)}$ $H$ is measurable function,
$\int_{\mathbb{R}^d}H(z)p(z)dz\equiv 0$; for $\alpha\ge 1$, sufficiently
small $\delta>0$ and any $\beta<0\wedge\frac{1}{2}(3-\alpha-\delta)$,
\[
\|H(x)\|\le \mathfrak{c}(1+\|x\|)^{\alpha-1+\beta}.
\]
\begin{remark}\label{rem-1new}
{\rm Under $\mathbf{(A_b)}$, $\mathbf{(A_{\sigma, a})}$ and
$\mathbf{(A_H)})$, from Pardoux and Veretennikov, \cite{PVI} Theorem 2,
it follows that $U(x)$, given in \eqref{1.5U}, is bounded and
solves the Poisson equation
\begin{equation*}
\mathscr{L}U=-H
\end{equation*}
in the class of functions with Sobolev's partial second
derivatives locally integrable in any power and a polynomial
growth. With all this going on,
\begin{equation}\label{nabla2.3}
|\nabla U(x)|\le c(1+\|x\|)^{(\beta+\alpha-1)^+}
\end{equation}
and, by embedding theorems \cite{LSU}, all entries of $\nabla U$
are continuous functions. So, the Krylov generalization of It\^o's
formula (see \cite{KI}) is applicable to $U(X_t)$:
\begin{equation}\label{U}
U(X_t)=U(x)-\int_0^tH(X_s)ds+\int_0^t\nabla
U(X_s)\sigma(X_s)dW_s.
\end{equation}
}
\end{remark}
\begin{theorem}\label{theo-V1}
Under $\mathbf{(A_b)}$, $\mathbf{(A_{\sigma,a})}$ and
$\mathbf{(A_H)}$, the family $(S^\kappa_t)_{t\ge 0}$ obeys the MDP
in $(\mathbb{R}^q,\rho)$ with the rate of speed $\varrho(t)$ and
the rate function given in \eqref{1.7e} with $Q$ defined in
\eqref{QQQQQ}.
\end{theorem}

\subsection{Nonlinear model, II}
\label{sec-03}

Though Theorem \ref{theo-V1} serves a wide class of bounded and
unbounded functions $H$, it is far from to be universal especially
for $\alpha=1$.

So, we fix the next set of stronger assumptions.

\medskip
$\mathbf{(A'_{b,\sigma})}$ \ \ $b(x)$ and $\sigma(x)$ are
Lipschitz continuous; for any $x',x''\in \mathbb{R}^d$ there

\hskip .48in exists a positive number $\nu$ such that
\begin{gather*}
2\lef(x'-x'',b(x')-b(x'')\rig+\trace[\sigma(x')-\sigma(x'')][\sigma(x')-\sigma(x'')]^*
\\
\le -\nu\|x'-x''\|^2.
\end{gather*}

\medskip
$\mathbf{(A'_{a})}$ $ \lambda \mathrm{I}\leq a(z)\leq \Lambda
\mathrm{I}, $ for some $\Lambda>\lambda>0$.

\medskip
$\mathbf{(A'_H)}$ \hskip .04in  $H(x)$ is Lipschitz continuous
function.

\begin{theorem}\label{theo-new3.1}
Under $\mathbf{(A'_{b,\sigma})}$, $\mathbf{(A'_{a})}$ and
$\mathbf{(A'_H)}$, the statement of Theorem {\rm \ref{theo-V1}} remains
valid.
\end{theorem}

\subsection{Linear model}

The diffusion process $X_t$ solves \eqref{itolin}, $A=A_{d\times
d}$, $B=B_{d\times d}$ and $(W_t)_{t\ge 0}$ is a standard
vector-valued Wiener process of the corresponding size.

For this setting, $\mathbf{(A_b)}$
or $\mathbf{(A'_{b,\sigma})}$, and $\mathbf{(A_{a})}$ are too restrictive.
We replace them by the following assumptions.

\medskip
$\mathbf{(A)}$ Eigenvalues of $A$ have negative real parts.

\medskip
$\mathbf{(A_B)}$ $ D:=BB^*+A^*BB^*A+\ldots+(A^*)^{d-1}BB^*A^{d-1}
$ is nonsingular ma-

\hskip .4in trix.

\medskip
$\mathbf{(A''_{H})}$ Suppose either

\hskip .4in 1) $H$ possesses continuous and bounded partial
derivatives,

\hskip .4in 2) $H$ is bounded, $BB^*>0$.

\begin{theorem}\label{theo-2}
Under $\mathbf{(A)}$, $\mathbf{(A_B)}$ and $\mathbf{(A''_H)}$, the family $(S^\kappa_t)_{t\to \infty}$
obeys the MDP in $(\mathbb{R}^d,\rho )$ with rate of speed
$\varrho(t)$ and the rate function given in \eqref{1.7e} with $Q$
defined in \eqref{QQQQQ}.
\end{theorem}

The next result deals with quadratic function $H$. Under
$\mathbf{(A)}$ and $\mathbf{(A_B)}$, the invariant measure $\mu$ is zero
mean Gaussian with nonsingular covariance matrix $P$ solving the
Lyapunov equation

\begin{equation}\label{Lyapunov}
A^*P+PA+BB^*=0.
\end{equation}
We introduce also a positive definite matrix $\Gamma=\Gamma_{q\times q}$ and
a matrix $\Upsilon=\Upsilon_{q\times q}$ solving  the Lyapunov equation
\begin{equation*}
A^*\Upsilon+A\Upsilon+\Gamma=0.
\end{equation*}

\begin{theorem}\label{theo-3.4}
Assume $\mathbf{(A)}$ and $BB^*>0$ and
$$
H(x)=\lef x,\Gamma x\rig-\trace(\Gamma^{1/2}P\Gamma^{1/2}).
$$

Then, the family $(S^\kappa_t)_{t\to \infty}$ obeys the MDP in
$(\mathbb{R}^d,\rho )$ with rate of speed $\varrho(t)$ and the
rate function given in \eqref{1.7e} with
$$
Q=4\trace(\Upsilon BPB^*\Upsilon)>0.
$$
\end{theorem}

\subsection{More examples}

In this section, we give examples which are not explicitly
compatible with Theorems \ref{theo-V1} - \ref{theo-3.4}.

\begin{example}\label{ex-1}

 Let $d=1$, $H(x)=x^3$ and
\begin{equation}\label{3.1a}
dX_t=-X^3_tdt+dW_t.
\end{equation}
Though $\mathbf{(A_b)}$ holds with $\alpha=3$, Theorem \ref{theo-V1} is not
applicable since by $\mathbf{(A_H)}$ only $H$ with property $\|H(x)\|\le
c(1+\|x\|)^\gamma$, $\gamma<2$ is admissible.

Nevertheless, the MDP holds and is trivially verified. Indeed,
\eqref{3.1a} is nothing but \eqref{U} with $U(x)\equiv x$. Hence,
$\nabla U(x)=1$ and $Q=1$.

Consequently, (ii) from Theorem \ref{theo-main} automatically holds.

(i) is reduced to
$
\lim_{t\to\infty}\varrho(t)\log P\big(X^2_t\ge
t^{2\kappa}\varepsilon\big) =-\infty
$
and is verified with the help of Theorem \ref{theo-a.1}  with
$V(x)\equiv x^2$. Actually, by It\^o's formula we find that $
dV(X_t)=[-2V^2(X_t)+1]dt+dN_t, $ where $N_t=\int_0^t2X_sdW_s$.
Hence,
$$
\mathscr{L}V(x)\le -V^2(x)+1\quad\text{and}\quad  \langle
N\rangle_t=\int_0^t4V(X_s)ds.
$$
\end{example}

\begin{example}\label{ex-2}

Let $d=1$ and
$$
dX_t=b(x_t)dt+dW_t,
$$
where $b(x)$ is Lipschitz continuous and symmetric,
$
b(x)=-b(-x),
$
function (obviously $b(0)=0$).
Under $\mathbf{(A'_{b,\sigma})}$, providing $\mathbf{(A_b)}$,
$X_t$ is an ergodic diffusion process with the symmetric invariant density,
$p(z)=p(-z)$. So,  any bounded $H(x)$, with $H(x)=-H(-x)$, possesses
\eqref{nul}. We choose
$$
H(x)=\sign(x), \ \text{letting $\sign(0)=0.$}
$$
However, neither Theorem \ref{theo-2}  nor Theorem \ref{theo-V1} are compatible
with the setting owing to  $H(x)$ does not satisfy neither $\mathbf{(A'_H)}$ nor
$\mathbf{(A_H)}$.
Nevertheless, we show that the standard MDP holds.
A computational trick proposes to use a decomposition
$H=H'+H''$ for
$$
H'(x)=
  \begin{cases}
    e^{-x}, & x> 0
    \\
    0, & x=0
    \\
    -e^x, & x<0
  \end{cases}
$$
since $H'$ satisfies $\mathbf{(A_H)}$ and
$$
H''(x)=
  \begin{cases}
    1-e^{-x}, & x\ge 0
    \\
    -1+e^x, & x<0
  \end{cases}
$$
satisfies $\mathbf{(A'_H)}$. Then, $U'(x)$ and $\nabla U'(x)$ are well
defined and both are bounded; at the same time
$U''(x)$ and $\nabla U''(x)$ are also well defined and $\nabla U''(x)$ is bounded,
i.e.  $|U''(x)|\le c(1+|x|)$.

Taking $U(x)=U'(x)+U''(x)$ we get bounded $\nabla U(x)=\nabla U'(x)+\nabla U''(x)$ and
$U(x)$ satisfying the linear growth condition. Moreover, due to
$
M_t=M'_t+M''_t,
$
we have
$$
M_t=\int_0^t\nabla U'(X_s)dW_s+\int_0^t\nabla U''(X_s)dW_s=\int_0^t\nabla U(X_s)dW_s,
$$
providing $\langle M\rangle_t=\int_0^t\big(\nabla U(X_s)\big)^2ds$
with bounded $\big(\nabla U(x)\big)^2$.

Now, (i) and (ii) from Theorem \ref{theo-main} are verified in a
standard way with the help of Theorems \ref{theo-a.1}, \ref{theo-a.2}.
\end{example}

\begin{example}\label{ex-3}
({\sf Linear version of Langevin model.}) A nonlinear Langevin's
model, including our linear one, is studied in Wu, \cite{WuH}. The
result from \cite{WuH} seems not to be accomplished. At least, we
could not adapt assumptions from there to verify the MDP for the
following setting.

Let $X_t=\begin{pmatrix}
  q_t \\
  p_t \\
\end{pmatrix}\in\mathbb{R}^{2d}$ with ($q_t,p_t\in \mathbb{R}^d$) and
\begin{equation}\label{5.1a}
\begin{aligned}
dq_t&=p_tdt
\\
dp_t&=-\Gamma p_tdt- \nabla \mathrm{F}(q_t)dt+\sigma dW_t,
\end{aligned}
\end{equation}
where $\nabla \mathrm{F}(q)=\Lambda q$ and matrices $\Lambda$,
$\Gamma$ and $\sigma\sigma^*$ are positive definite. We verify the
MDP with the help of Theorem \ref{theo-2}.

It is expedient to write \eqref{5.1a} to the form \eqref{itolin}
with matrices (in a block form)
\begin{equation*}\label{AB}
A=\begin{pmatrix}
  0 & \mathrm{I} \\
  -\Lambda & -\Gamma \\
\end{pmatrix}
\quad\text{and}\quad B=\begin{pmatrix}
  0 & 0 \\
  0 & \sigma \\
\end{pmatrix}.
\end{equation*}

\noindent In accordance with Theorem \ref{theo-2}, we have to
verify only two conditions:

1) eigenvalues of $A$ have negative real parts,

2) the matrix $\mathrm{D}$ (see $\mathbf{(A_B)}$) is nonsingular.

\medskip
\noindent 1) fulfils since free of noise \eqref{5.1a}:
\begin{equation*}
\begin{aligned}
\dot{q}_t&=p_t
\\
\dot{p_t}&=-\Gamma p_t- \nabla \mathrm{F}(q_t)
\end{aligned}
\end{equation*}
is asymptotically stable. Traditionally for the Langevin equation,
this result is easily verified with the help of Lyapunov's function
$V_t=\frac{1}{2}\|p_t\|^2+F(q_t)$ and is omitted here.

2) holds since $ D':=BB^*+A^*BB^*A(\le D) $ is nonsingular.
Indeed,
$$
D'=\begin{pmatrix}
  \Lambda\sigma^*\sigma\Lambda & \Lambda\sigma^*\sigma\Gamma \\
  \Gamma\sigma^*\sigma\Lambda & \Gamma\sigma^*\sigma\Gamma+\sigma^*\sigma
\end{pmatrix},
$$
that is, with a vector $ v=\begin{pmatrix}
  v_1 \\
  v_2 \\
\end{pmatrix}\ne 0,
$ we have
$$
\begin{aligned}
\lef v,D^*Dv\rig&=\lef v_1,\Lambda\sigma^*\sigma\Lambda v_1\rig +
\lef v_2,\Gamma\sigma^*\sigma\Gamma v_2\rig+ 2\lef
v_1,\Lambda\sigma^*\sigma\Gamma v_2\rig
\\
&\quad +\lef v_2,\sigma^*\sigma v_2\rig.
\end{aligned}
$$
By virtue of the well known inequality
\begin{equation}\label{n}
2\lef z_1,z_2\rig\ge -\lef z_1,z_1\rig-\lef z_2,z_2\rig,
\end{equation}
we get $ \lef v_1,\Lambda\sigma^*\sigma\Lambda v_1\rig+ \lef
v_2,\Gamma\sigma^*\sigma\Gamma v_2\rig+ 2\lef
v_1,\Lambda\sigma^*\sigma\Gamma v_2\rig\ge 0. $ Consequently,
under $v_2\ne 0$, we have $ \lef v,(D')^*D'v\rig \ge \lef
v_2,\sigma^*\sigma v_2\rig>0. $ Even though $v_2=0$, and so
$v_1\ne 0$, we also have $ \lef v,(D')^*D'v\rig= \lef
v_1,\Lambda\sigma^*\sigma v_1\Lambda\rig>0$.

Thus, under $\mathbf{(A''_H)}$,  the MDP holds.
\end{example}

\begin{example}\label{ex-4}
({\sf MDP for a smooth component of diffusion process}.) Let
$X^{(1)}_t$ be the first component of a diffusion process $X_t$ with entries
$X^{(i)}_t$, $i=1,\ldots,d$:
\begin{equation}\label{3.5nn}
\begin{aligned}
\dot{X}^{(1)}_t&=X^{(2)}_t
\\
\dot{X}^{(i)}_t&=X^{(i+1)}_t, \ i=2,\ldots,d-1
\\
dX^{(d)}_t&=-\sum_{i=1}^da_iX^{(d-i)}_tdt+bdW_t,
\end{aligned}
\end{equation}
where $a_1,a_2,\ldots,a_d$ and $b$ are positive numbers and $W_t$
is a Wiener process.

As in  the previous example, we rewrite
\eqref{3.5nn} to the form of \eqref{itolin} with
$$
A=
\begin{pmatrix}
  A_{11} & A_{12} \\
  A_{21} & A_{22} \\
\end{pmatrix}
\quad\text{and}\quad B=
\begin{pmatrix}
  B_{11} & B_{12} \\
  B_{21} & B_{22} \\
\end{pmatrix},
$$
where
$$
A_{11}=
\begin{pmatrix}
  0 & 1 & 0 &\ldots &0 &0 \\
  0 & 0 & 1 & 0 & \ldots &0 \\
  \ldots & \ldots & \ldots & \ldots & \ldots & \ldots \\
  0 & 0 & 0 & 0 & 0 & 1
\end{pmatrix}_{(d-1)\times(d-1)},
$$
$A_{12}=0$, $A_{22}=-a_d$,
$
A_{21}=\begin{pmatrix}
  -a_1 & -a_2 & \ldots & -a_{d-1}
\end{pmatrix}_{1\times(d-1)}
$
and, analogously, $B_{11}=0_{(d-1)\times(d-1)}$, $B_{12}=0$, $B_{22}=b$,
$B_{21}=0_{1\times(d-1)}$.

We verify the
MDP with the help of Theorem \ref{theo-2}.
In order to guarantee $\mathbf{(A_H)}$, suffice it to assume that
roots of  the polynomial
$$
\phi(z)=z^d+a_1z^{d-1}+\ldots+a_{d-1}z+a_{d}
$$
have negative real parts owing to the noise free version of \eqref{3.5nn} is nothing but the
differential equation
$
x^{(d)}_t+\sum_{i=1}^{d-1}a_ix^{(d-i)}_t+a_dx_t=0.
$

Notice that
$\mathbf{(A_B)}$ is fulfilled too since
$D'=BB^*+A^*BB^*A (\le D)$ is a
nonsingular matrix. Actually,
$
D'=b^2\begin{pmatrix}
  A^*_{21}A_{21} & A^*_{21}A_{22} \\
  A_{21}A_{22} & A^2_{22}+1 \\
\end{pmatrix}
$
and so, we have
\begin{equation*}
\begin{aligned}
\lef v,D^*Dv\rig=b^2\Big[v^2_1\|A_{21}\|^2 +(A^2_{22}+1)\|v_2\|^2
+2v_1A_{22}\lef v_2,A_{21}\rig\Big].
\end{aligned}
\end{equation*}
Taking $ v=\begin{pmatrix}
  v_1 \\
  v_2 \\
\end{pmatrix}\ne 0,
$ where $v_1$ is a number and $v_2$ is a vector of the size $d-1$,
for $v_2=0$, and then $v_1\ne 0$, we have $\lef v,D^*Dv\rig>0$. Even though
$v_2\ne 0$, the use of $\eqref{n}$ provides $ \lef v,D^*Dv\rig\ge
b^2A^2_{22}\|v_2\|^2>0.$

\medskip
In order to establish the MDP for the family
$
\big(\frac{1}{t^\kappa}\int_0^tH(X^{(1)}_s)ds\big)_{t\to\infty},
$
we redefine the function $H$ as:
$
H(x^{(1)})\equiv \mathsf{H}(x^{(1)},x^{(2)},\ldots,x^{(d)})
$
and assume $\mathsf{H}$ satisfies $\mathbf{(A''_H)}$. Then, the family
$\big(\frac{1}{t^\kappa}\int_0^t\mathsf{H}(X_s)ds\big)_{t\to\infty}$
obeys the MDP with the rate of speed  $\varrho(t)$ and the rate
function
$
\mathsf{J}(Y)=\mathsf{J}(Y^{(1)},\ldots,Y^{(d)})
$
of the standard form \eqref{1.7e}.

Now, the desired MDP holds by Varadhan's contraction principle,
\cite{V}, with the same rate of speed and the rate
function
$$
j(y)=\inf\limits_{\{Y^{(2)},\ldots,Y^{(d)}\}\in\mathbb{R}^{d-1}}
J(y,Y^{(2)},\ldots,Y^{(d)}).
$$
\end{example}

\begin{example}\label{ex-5}
Let $X_t \ (\in\mathbb{R})$ be Gaussian diffusion with
\begin{equation*}
dX_t=-X_tdt+dW_t
\end{equation*}
and $H(x)=x^2\sign(x)$. This function satisfies \eqref{nul} and, at the same time,
is not compatible
with Theorems \ref{theo-V1} - \ref{theo-3.4}. So, we suppose to embed this setting to
a new one with a vector function $\mathsf{H}(x)$ with entries:
$$
\mathsf{H}_1(x)=\frac{1}{2}\sign(x)\quad\text{and}\quad
\mathsf{H}_2(x)=x^2\sign(x)-\frac{1}{2}\sign(x),
$$
which is MDP verifiable.
Applying arguments from the proof of the Theorem \ref{theo-2}, one can show
the existence of
$U_1(x)$ with bounded $\nabla U_1(x)$ such that
\begin{equation*}
\begin{aligned}
U_1(X_t)&=U_1(x)-\int_0^t\mathsf{H}_1(X_s)ds+M^{(1)}_t
\\
\langle M^{(1)}\rangle_t&=\int_0^t(\nabla U_1(X_s))^2ds.
\end{aligned}
\end{equation*}
Now, we establish similar property of $H_2(x)$.
By the Krylov-It\^o formula (see \cite{KI}), we find
that
$$
dH(X_t)=-\mathsf{H}_2(X_t)dt+|X_t|dW_t.
$$
Consequently, $U_2(x)\equiv H(x)$ and $\langle M^{(2)}\rangle_t=\int_0^2X^2_sds$.

Now, we may verify (i), (ii) from Theorem \ref{theo-main}.

(i):  Since $\nabla U_1$ is bounded, $U_1$ satisfies the linear growth condition. Thus,
$$
|U_1|\le c(1+|U_2|)=c(1+|H(x)|)\le \mathfrak{c}(1+x^2).
$$
Hence, (i) is reduced to
$
\lim_{t\to\infty}\varrho(t)\log P\big(X^2_t\ge
t^\kappa \varepsilon\big)=-\infty.
$
The latter holds owing to $X^2_t$ possesses an exponential
moment: $Ee^{\lambda X^2_t}<\infty$  uniformly
in $t$ over $\mathbb{R}_+$ and sufficiently small $\lambda$ and, therefore,
the Chernoff inequality is effective. Write
\begin{equation*}
\begin{aligned}
\frac{1}{t^{2\kappa-1}}\log P\big(X^2_t>t^\kappa\varepsilon\big)&\le
\frac{1}{t^{2\kappa-1}}\log \Big(e^{-\lambda
t^\kappa\varepsilon+ \log Ee^{\lambda X^2_t}}\Big)
\\
&\le -\lambda t^{1-\kappa}\varepsilon+\frac{\log Ee^{\lambda
X^2_t}}{t^{2\kappa-1}} \xrightarrow[t\to\infty]{} -\infty.
\end{aligned}
\end{equation*}

Notice that $|U_2(x)|=x^2$, so that, the (i) verification is the same as for $U_1$.

\medskip
(ii): The martingale $M_t$ is vector-valued process with two entries $M^{(1)}_t$ and
$M^{(2)}_t$. Hence, its variation process is a matrix
$$
\langle M\rangle_t=\begin{pmatrix}
  \langle M^{(1)}\rangle_t & \langle M^{(1)},M^{(2)}\rangle_t \\
\langle M^{(1)},M^{(2)}\rangle_t &   \langle M^{(2)}\rangle_t,
\end{pmatrix}
$$
so that, the entries of $Q$ are defined in the following way:
\begin{equation*}
\begin{aligned}
Q_{11}&=\int_\mathbb{R}\big(\nabla U_1(z)\big)^2p(z)dz,
\\
Q_{22}&=\int_\mathbb{R}z^2p(z)dz,
\\
Q_{12}&=\int_\mathbb{R}\nabla U_1(z)|z|p(z)dz.
\end{aligned}
\end{equation*}

Thus, (ii) is reduced to
\begin{equation}\label{3.10ee}
 \lim_{t\to\infty}\varrho(t)\log
P\Big(\Big|\int_0^th(X_s)ds\Big|\ge
t\varepsilon\Big)=-\infty,
\end{equation}
where $h(x)$ is continuous function satisfying \eqref{nul} and
either

1) $|h(x)|\le c|x|$,

2) $h(x)=x^2-\frac{1}{2}$.

\smallskip
\noindent
In 1), we apply
$
h(x)=h'_l(x)+h''_l(x),
$
borrowed from \eqref{2.7m}, and verify versions of \eqref{3.10ee} with
$h'_l$ and $h''_l$ separately.

$h'_l$-version holds owing to by Theorem \ref{theo-V1} ($\alpha =1$)
$
\Big(\frac{1}{t^\kappa}\int_0^th'_l(X_s)ds\Big)_{t\to\infty}
$
obeys $\varrho$-MDP with a nondegenerate rate function and $\kappa<1$.

$h''_l$-version holds owing to
$
|h''_l(x)|\le I(|x|>l)|x|\le \frac{x^2}{l}
$
and for sufficiently large $l$,
$
\lim_{t\to\infty}\varrho(t)\log P\Big(\int_0^tX^2_sds>tl\varepsilon\Big)
=-\infty
$
verified with the help of Theorem \ref{theo-a.1} for $V(x)=x^2$.

In 2), by Theorem \ref{theo-3.4}, $\big(\frac{1}{t^\kappa}
\int_0^t\big[X^2_s-\frac{1}{2}\big]ds\big)_{t\to\infty}$ obeys $\varrho$-MDP
with a nondegenerate rate function. So, it remains to recall that
$\kappa<1$.

\medskip

Thus, $\varrho$-MDP for new family holds true with the rate function
$J(Y)$, $Y\in \mathbb{R}^2$, defined in \eqref{1.7e}. Hence, the original
family possesses the MDP with the quadratic rate function
$$
j(y)=\inf_{\{Y_1,Y_2:Y_1+Y_2=y\}}J(Y).
$$
\end{example}

\section{Proof of Theorems from Section \ref{sec-3new}}
\label{sec-4}

\medskip
\noindent

\subsection{The proof of Theorem \ref{theo-V1}}
\label{sec-2.2}

Denote by
$
M_t=\int_0^t\nabla U(X_s)\sigma(X_s)dW_s
$
the martingale from \eqref{U} having
$ \langle M\rangle_t=\int_0^t\nabla U(X_s)a(X_s)\nabla^*U(X_s)ds$.

We shall verify (i) and (ii) from Theorem \ref{theo-main}.

\medskip
(i) holds since, by Remark \ref{rem-1new}, $U$ is bounded.

\medskip
(ii) is verified in a few steps.

\noindent
\underline{Step 1: $Q$ identification. }
We show that $\int_{\mathbb{R}^d}\nabla U(z)a(z)\nabla^*U(z)p(z)dz=Q.
$
This fact is well known and is given here for a reader convenience only.
Notice that, by \eqref{1.5U},
\begin{equation*}
Q=E\big[H(X^\mu_0)U^*(X^\mu_0)+U(X^\mu_0)H^*(X^\mu_0)\big],
\end{equation*}
where $X^\mu_t$ the stationary version of $X_t$, that is, the version
solving \eqref{ito00} subject to $X^\mu_0$ the random
vector, independent of $W_t$, with the distribution
provided by the invariant measure $\mu$. Hence, suffice it to show that
\begin{equation}\label{x3}
E\big[\nabla U(X^\mu_0)a(X^\mu_0)\nabla^*U(X^)\mu_0\big]=
E\big[H(X^\mu_0)U^*(X^\mu_0)+U(X^\mu_0)H^*(X^\mu_0)\big].
\end{equation}
We verify \eqref{x3} with the help of It\^o's formula
\begin{gather*}
U(X^\mu_t)U^*(X^\mu_t)=U(X^\mu_0)U^*(X^\mu_0)
-\int_0^t\big[H(X^\mu_s)
U^*(X^\mu_s)+U(X^\mu_s)H^*(X^\mu_0)\big]ds
\\
+\int_0^t\big[U(X^\mu_s)dM^*_s+dM_sU^*(X^\mu_0)\big]
+\int_0^t\nabla(X^\mu_s)a(X^\nu_s)\nabla^*(X^\mu_s)ds
\end{gather*}
by taking the expectation.

\medskip
\underline{Step 2. Preliminaries.} Set
$ \mathsf{H}(x)=\nabla U(x)a(x)\nabla^*(x)-Q$ and let $h(x)$ denotes
any entry of $\mathsf{H}(x)$. For (ii) to be valid suffice
it to show that
\begin{equation}\label{to0}
\lim_{t\to\infty}\varrho(t)\log
P\Big(\Big|\int_0^th(X_s)ds\Big|>t\varepsilon\Big)=-\infty.
\end{equation}
Recall $\int_{\mathbb{R}^d}h(z)p(z)dz=0$. By \eqref{nabla2.3},
\begin{equation*}
|h(x)|\le c(1+\|x\|)^{2(\beta+\alpha-1)^+}.
\end{equation*}

\medskip
\noindent We consider separately two  cases provided by a special
choice
of
$$
\beta<0\wedge \frac{1}{2}(3-\alpha-\delta) \ \text{for sufficiently small
$\delta>0$}.
$$
(see, $\mathbf{(A_H)}$):

- $(\alpha=1):$ $|h(x)|$ is bounded;

- $(\alpha>1):$ $|h(x)|\le c(1+\|x\|)^{1+\alpha-\delta}$, $1+\alpha-\delta\ge 2.$

\medskip
\underline{Step 3.} \fbox{$\alpha=1$} For sufficiently large number $l$, set
\begin{equation*}\label{2.7m}
h'_l(x)=
  \begin{cases}
    h(x) & \|x\|\le l,
    \\
    v_l(x) & l<\|x\|\le l+1
    \\
    0 & \|x\|>l+1,
  \end{cases}
  \end{equation*}
where $v_l(x)$ is bounded continuous function such that $h'_l(x)$
is continuous function with $\int_{\mathbb{R}^d}h'_l(z)p(z)dz=0$.
In contrast to $h$, the function $h'_l$  decreases fast  to zero
with $\|x\|\to\infty$, so that, a negative constant $\beta'$ can
be chosen such that
$$
|h'_l(x)|\le c(1+|x|)^{\beta'+\alpha-1}\equiv c(1+|x|)^{\beta'}.
$$
In accordance with this property, $ u(x)=-\int_0^\infty
E_xh'(X_t)dt $ solves the Poisson equation $ \mathscr{L}u=-h'_l $
and is bounded jointly with $\nabla u(x)$ (see, Remark
\ref{rem-1new}). Hence, $ u(X_t)=u(x)-\int_0^th'(X_s)ds+m_t $ with
the martingale $m_t=\int_0^t\nabla u(X_s)\sigma(X_s)dW_s$ having
$
\langle m\rangle_t=\int_0^t\nabla u(X_s)a(X_s) \nabla^* u(X_s)ds.
$
The negligibility of $\frac{u(x)-u(X_t)}{t}$ in $\varrho$-MDP scale is
provided by the boundedness of $u(x)$. The same type
negligibility of $\frac{1}{t}m_t$ is provided by
the boundedness of $\nabla u^*(x)a(x)\nabla u(x)$, due to Theorem \ref{theo-a.2}.

Consequently, a version of \eqref{to0} with $h'_l$ holds true.

\medskip
Set $h''_l=h-h'_l$. Since $h$ is bounded, $|h''_l(x)|\le
cI(\|x\|>l)\le \frac{c}{l^2}\|x\|^2.$
Consequently a version of \eqref{to0} with $h''_l$ is reduced
to
\begin{equation*}
\lim_{t\to\infty}\varrho(t)\log P\Big(
\int_0^t\|X_s\|^2ds>t(l^2\varepsilon)\Big)=-\infty
\end{equation*}
and is verified with the help of Theorem \ref{theo-a.2} for $V(x)=\|x\|^2$
owing to
$$
\mathscr{L}V(x)\le -cV(x)+\mathfrak{c},
\quad\text{and}\quad
\langle N_t\rangle\le\int_0^t\mathbf{c}\big(1+V(X_s)\big)ds
$$
are fulfilled under $\mathbf{(A_b)}$ and $\mathbf{(A_{\sigma, a})}$ (a verification
of these facts is accomplished with the help of It\^o's formula).

\medskip
\underline{Step 4.} \fbox{$\alpha>1$}
We apply again the decomposition $h=h_l'+h_l''$. With chosen $l$,
$|h'_l|$ is decreasing fast to zero, with $\|x\|\to\infty$, and is bounded by
$c(1+l)^{1+\alpha-\delta}$. So, the version of \eqref{to0} with $h'_l$
is verified as in the case ``$\alpha=1$''.

\medskip
Notice that
$$
|h_l''(x)|\le c(1+\|x\|)^{1+\alpha-\delta)}I(\|x\|>l)\le
\frac{c}{l^\delta} (1+\|x\|)^{1+\alpha}\le\frac{\mathfrak{c}}{l^\delta}\big(1+V(x)\big),
$$
where
$
V(x)=\frac{\|x\|^{4+2\alpha}}{1+\|x\|^{3+\alpha}}.
$
Hence, the version of \eqref{to0} with $h''_l$ is reduced to
\begin{equation*}
\lim_{t\to\infty}\varrho(t)\log
P\Big(\int_0^tV(X_s)ds>tl^\delta\Big)=-\infty.
\end{equation*}
To this end, we apply Theorem \ref{theo-a.1}.

First, taking into account that
$\|x\|^{3+\alpha}=\big(\|x\|^2\big) ^{\frac{3+\alpha}{2}}$,
$\|x\|^{4+2\alpha}=\big(\|x\|^2)^{2+\alpha}$ and
$\frac{3+\alpha}{2}>2$, by the It\^o formula we find that
\begin{equation*}
\begin{aligned}
d\|X_t\|^2&=\big[2\lef X_t,b(X_t\rig+\trace\big(a(X_t)\big)\big]dt+
2\lef X_t,\sigma(X_t)dW_t\rig,
\\
d\|X_t\|^{3+\alpha}&=\Big(\frac{3+\alpha}{2}-1\Big)\big(\|X_t\|^2\big)^{\frac{3+\alpha}{2}-1}
\big\{2\lef X_t,b(X_t)\rig +\trace\big(a(X_t)\big)\big\}
\\
&\quad
+2\Big[\frac{3+\alpha}{2}-1\Big]\Big[\frac{3+\alpha}{2}-2\Big]
\big(\|X_t\|^2\big)^{\frac{3+\alpha}{2}-2}\lef
X_t,a(X_t)X_t\rig\Big]dt
\\
&\quad
+\Big(\frac{3+\alpha}{2}-1\Big)\big(\|X_t\|^2)^{\frac{3+\alpha}{2}-1}
2\lef X_t,\sigma(X_t)dW_t\rig,
\\
d\|X_t\|^{4+2\alpha}&=\big(1+\alpha\big)\big(\|X_t\|^2\big)^{1+\alpha}
\big\{2\lef X_t,b(X_t)\rig +\trace\big(a(X_t)\big)\big\}
\\
&\quad +2\big[1+\alpha\big]\alpha \big(\|X_t\|^2\big)^\alpha\lef
X_t,a(X_t)X_t\rig\Big]dt
\\
&\quad +\big(1+\alpha\big)\big(\|X_t\|^2)^{1+\alpha} 2\lef
X_t,\sigma(X_t)dW_t\rig,
\\
d\frac{1}{1+\|X_t\|^{3+\alpha}}&=-\frac{d\|X_t\|^{3+\alpha}}{(1+\|X_t\|^{3+\alpha})^2}
\\
&\quad+\frac{2(1+\alpha)\|X_t\|^{1+\alpha}}{(1+\|X_t\|^{3+\alpha})^3}\lef
X_t, a(X_t)X_t\rig dt,
\\
dV(X_t)&=\frac{d\|X_t\|^{4+2\alpha}}{1+\|X_t\|^{3+\alpha}}+\|X_t\|^{4+2\alpha}
d\frac{1}{1+\|X_t\|^{3+\alpha}}
\\
&\quad +\frac{2\big(1+\alpha\big)^2\|X_t\|^{3(1+\alpha)}\lef
X_t,a(X_t)X_t\rig} {(1+\|X_t\|^{3+\alpha})^2}dt.
\end{aligned}
\end{equation*}

Thus, we have $dV(X_t)=\mathscr{L}V(X_t)dt+dN_t$, where
\begin{equation*}
\begin{aligned}
\mathscr{L}V(x)&=\frac{1}{1+\|x\|^{3+\alpha}}\Big[(1+\alpha)\|x\|^{2(1+\alpha)}
\big\{2\lef x,b(x)x\rig+\trace(a(x))\big\}
\\
&\quad +2\alpha\big(1+\alpha)\|x\|^{2\alpha}\lef x,a(x)x\rig\Big]
\\
&\quad
-\frac{\|x\|^{4+2\alpha}}{{(1+\|X_t\|^{3+\alpha})^2}}\Big[\frac{1}{2}(1+\alpha)\|x\|^{1+\alpha}\{2\lef
x,b(x)x\rig+\trace(a(x))\}
\\
&\quad +\frac{1}{2}(1+\alpha)(\alpha-1)\|x\|^{\alpha-1}\lef
x,a(x)x\rig\Big]
\\
&\quad
+\frac{\|x\|^{4+2\alpha}}{{(1+\|X_t\|^{3+\alpha})^3}}2(1+\alpha)\|x\|^{1+\alpha}
\lef x,a(x)x\rig
\\
& \le (1+\alpha)\lef
x,b(x)x\rig\Big[\frac{2\|x\|^{2(1+\alpha)}}{1+\|x\|
^{3+\alpha}}-\frac{\|x\|^{4+2\alpha}}{(1+\|x\|^{3+\alpha})^2}\Big]+o(\|x\|^{2\alpha})
\\
&=(1+\alpha)\lef
x,b(x)x\rig\frac{2\|x\|^{2(1+\alpha)}+2\|x\|^{5+3\alpha}-
\|x\|^{4+2\alpha}}{(1+\|x\|^{3+\alpha})^2}+o(\|x\|^{2\alpha})
\\
&\le -c\|x\|^{2\alpha}+\mathfrak{c}\le
-\mathbf{c}V^{\frac{2\alpha}{1+\alpha}}(x)+ \mathfrak{c}
\end{aligned}
\end{equation*}
and
$
  N_t=\int_0^t\lef X_s,\sigma(X_s)dW_s\rig\Big[\frac{2(1+\alpha)\|X_s\|^{2(1+\alpha)}}
  {1+\|X_s\|^{3+\alpha}}
  -\frac{(1+\alpha)\|X_s\|^{5+3\alpha}}{(1+\|X_s\|^{3+\alpha})^2}\Big],
$
that is,
$$
\langle N\rangle_t\le
\int_0^t\big(c\|X_s\|^{2\alpha}+\mathfrak{c}\big)ds \le
\int_0^t\big(\mathbf{c}V^{\frac{2\alpha}{1+\alpha}}(X_s)+\mathfrak{c}\big)ds.
$$

Thus, the assumptions of Theorem \ref{theo-a.1} are fulfilled and,
thereby, for sufficiently large $l$, we have
$
\lim_{t\to\infty}\varrho(t)\log
P\big(\int_0^tV^{\frac{2\alpha}{1+\alpha}}(X_s)ds
>tl^\delta\big)=-\infty
$
and, it is left to notice that $\frac{2\alpha}{1+\alpha}>1$ for
$\alpha>1$. \qed

\subsection{The proof of Theorem \ref{theo-new3.1}}
\label{sec-4.2}

By $\mathbf{(A'_{b,\sigma})}$,
\begin{equation*}
\begin{aligned}
\lef x,b(x)\rig&=\lef x,(b(x)-b(0)\rig+\lef x,b(0)\rig
\\
& \le -\lef x,B_0x\rig+\|b(0)\|\|x\|
\\
& \le -\nu \|x\|^2+\|b(0)\|\|x\|
\end{aligned}
\end{equation*}
that is, there exists $r>0$ such that
$
\lef x,b(x)\rig\le -r\|x\|^2.
$

Hence,
$\mathbf{(A'_{b,\sigma})}\Rightarrow\mathbf{(A_{b})}(\alpha=1)$.
However since, by $\mathbf{(A'_H)}$,
$
\|H(z)\|\le c(1+\|z\|)
$
is admissible, Theorem 2 from
Pardoux and Veretennikov, \cite{PVI}, is no longer applicable.
At the same time, Theorem 1 from \cite{PVI} states that $U$ from \eqref{1.5U} solves
the Poisson equation $ \mathscr{L}U(z)=-H(z) $ and satisfies the
following properties: for some $m>2$,
$$
\|U(x)\|\le c(1+\|x\|^m)\quad\text{and}\quad \|\nabla U(x)\|\le
c(1+\|x\|^m).
$$
Nevertheless, regardless of that, $\mathbf{(A'_H)}$ provides
\begin{equation}\label{3.2new}
\|U(x)\|\le c(1+\|x\|)\quad\text{and}\quad \|\nabla U(x)\|\le c.
\end{equation}
Actually, let $X^x_t$ denotes the solution of \eqref{ito00} subject to
$X_0=x$. Since for any $x'$ and $x''$, we have  $ U(x')-U(x'')=\int_0^\infty
E\big[H(X^{x'}_t)-H(X^{x''}_t)\big]dt$, by
$\mathbf{(A'_H)}$, we have ($L$ is the Lipschitz constant for $H$)
\begin{equation*}
\begin{aligned}
|U(x')-U(x'')|&\le L\int_0^\infty
\big|E\big[X^{x'}_t-X^{x''}_t\big]\big|dt
\\
&\le L\int_0^\infty \Big(E\|X^{x'}_t-X^{x''}_t\|^2\Big)^{1/2}dt,
\end{aligned}
\end{equation*}
where
$
d[X^{x'}_t-X^{x''}_t]=[b(X^{x'}_t)-b(X^{x''}_t)]dt+[\sigma(X^{x'}_t)
-\sigma(X^{x''}_t)]dW_t.
$
With the help of It\^o's formula, we find that
\begin{equation*}
\begin{aligned}
d\|X^{x'}_t-X^{x''}_t\|^2_t&=2\lef
(X^{'}_t-X^{x''}_t),[b(X^{x'}_t)-b(X^{x''}_t)]\rig dt
\\
&\quad +2\lef
(X^{'}_t-X^{x''}_t),[\sigma(X^{x'}_t)-\sigma(X^{x''}_t)]dW_t\rig
\\
&\quad +\trace [\sigma(X^{x'}_t)-\sigma(X^{x''}_t)]
[\sigma(X^{x'}_t)-\sigma(X^{x''}_t)]^*.
\end{aligned}
\end{equation*}
Hence, $v_t=E\|X^{x'}_t-X^{x''}_t\|^2$ is differentiable relative
to $dt$ and
\begin{multline*}
\dot{v}_t=2E\Big[\lef
[X^{x'}_t-X^{x''}_t],[b(X^{x'}_t-b(X^{x''}_t)]\rig
\\
+\trace[\sigma(X^{x'}_t)-\sigma(X^{x''}_t)][\sigma(X^{x'}_t)-\sigma(X^{x''}_t)]^*
\Big].
\end{multline*}
Then, by $\mathbf{(A'_{b,\sigma})}$, we have $ \dot{v}_t\le -\nu
v_t, $ i.e.,
$
v_t\le \|x'-x''\|^2e^{-t\nu}.
$
The latter implies the Lipschitz continuity of $U$ and, in turn,
\eqref{3.2new}.

\medskip
We proceed the proof with the verification of (i) and (ii) from
Theorem \ref{theo-main}.

\medskip
(i): Due to \eqref{3.2new}, suffice it show that
$
\lim_{t\to\infty}\varrho(t)\log P\big(\|X_t\|^2>\varepsilon
t^{2\kappa} \big) =-\infty
$
what is verified with the help of Theorem \ref{theo-a.1} for $V(x)=\|x\|^2$.
With the help of It\^o's formula, one can find that
$$
\mathscr{L}V(x)=2\lef x,b(x)\rig+\trace a(x) \ \ \text{and} \ \
N_t=\int_0^t2\lef X_s,\sigma(X_s)dW_s\rig
$$
and next that $\mathscr{L}V(x)\le-cV(x)+\mathfrak{c}, \quad \langle
N\rangle_t\le\int_0^t\mathbf{c}V(X_s)ds$.

\medskip
(ii): It is verified similarly to \eqref{to0} for $\alpha=1$.

\qed

\subsection{The proof of Theorem \ref{theo-2}}
\label{sec-4.33}
Under $\mathbf{(A)}$, $\mathbf{(A_{B})}$, the
Pardoux-Veretennikov concept is no longer valid. Nevertheless,
$\mathbf{(A)}$ and  $\mathbf{(A_B)}$ provide the ergodicity of
$X=(X_t)_{t\ge 0}$ with the unique zero mean Gaussian invariant
measure characterized by a nonsingular covariance matrix $P$ solving
Lyapunov's equation, see \eqref{Lyapunov}.

We prove the theorem in a few steps.

\medskip
\underline{Step 1. Invariant and transition densities.} For $X_0=x$,
the diffusion process $X_t$ is Gaussian with the
expectation $EX_t=e^{At}x$ and the covariance matrix
$$
\cov(X_t,X_t)=\int_0^te^{(t-s)A^*}BB^*e^{(t-s)A}ds=:P_t
$$
solving the differential equation
\begin{equation}\label{Lyapt}
\dot{P}_t=A^*P_t+P_tA+BB^*
\end{equation}
subject to $P_0=0$.
It is well known, and is readily verified that, under
$\mathbf{(A)}$ and $\mathbf{(A_B)}$, we have
$
P_t>0
$
over $t>0$ and
$
\lim_{t\to\infty}P_t=P (>0).
$
If in addition $BB^*>0$, then, for $t$ in a vicinity of zero,
\begin{equation}\label{ebya}
|P^{-1/2}_t|\le \frac{c}{\sqrt{t}}.
\end{equation}
Since $P, P_t>0$, the invariant density $p(y)$ and the density of
$P^{(t)}_x(dy)$ relative to $dy$ are defined as:

\begin{gather*}
p(y)=\frac{1}{(2\pi \det P)^{d/2}}e^{-\frac{1}{2}\|y\|^2_{P^{-1}}}
\\
p(x,t,y)=\frac{1}{(2\pi \det
P_t)^{d/2}}\exp\Big(-\frac{1}{2}\big\|y-e^{tA}x
\big\|^2_{P^{-1}_t}\Big)
\end{gather*}

\medskip
\underline{Step 2. $U$ existence.} We prove that $U(x)$ from \eqref{1.5U} is well
defined
over $\mathbb{R}^d$ by showing
\begin{equation}\label{4.3<}
\int_0^\infty |E_xH(X_t)|dt<\infty.
\end{equation}

\medskip
Assume $\mathbf{(A''_H)}_{1)}$. Let  $X^\mu_t$, $X^x_t$ denote the stationary version of
$X_t$ and $X_t$ with $X_0=x$ respectively. By \eqref{nul} and the Lipschitz property
of $H$ (with the Lipschitz constant $L$), it holds
$
|E_x(X_t)|=|E[H(X^x_t)-H(X^\mu_t)|\le LE|X^x_t-X^\mu_t|,
$
where, by \eqref{ito00},
$
\frac{d}{dt}[X^x_t-X^\mu_t]=A[X^x_t-X^\mu_t],
$
i.e.,
$
[X^x_t-X^\mu_t]=e^{tA}[x-X^\mu_0].
$
Hence and by $\mathbf{(A)}$, there exists a positive constant $\lambda$ such that
$
|X^x_t-X^\mu_t|\le e^{-t\lambda}c(1+\|x\|+\|X^\mu_0\|).
$
The random vector $X^\mu_0$ is Gaussian, so that, $E\|X^\mu_0\|=c.$

Thus, $|E_x(X_t)|\le e^{-t\lambda}\mathfrak{c}(1+\|x\|)$ and
\eqref{4.3<} holds true.

\medskip
Assume $\mathbf{(A''_H)}_{2)}$.
We may adapt the results of Meyn and Tweedie, \cite{metw} (see also
Mattingly and Stuart, \cite{mast} and Mattingly Stuart and Higham, \cite{masth})
for getting \eqref{4.3<}.
However, taking into account the explicit formulae for $p(y)$ and
$p(x,t,y)$, the direct proof of \eqref{4.3<} is given.

For a definiteness, let $|H|\le K$. We apply an obvious inequality
$$
|E_xH(X_t)|\le K\int_{\mathbb{R}^d}\big|p(x,t,y)-p(y)\big|dy \
(\le 2K).
$$
A changing of variables: $z=(y-e^{tA}z)P^{-1/2}_t$ and the identity
\begin{equation*}
\begin{aligned}
& \frac{p(P^{1/2}_tz+e^{tA}x)}{p(x,t,P^{1/2}_tz+e^{tA}x)} =
\sqrt{\frac{\det P_t}{\det P}}
\\
&\times \exp\Big(-\frac{1}{2}\Big[\lef z,(P_tP^{-1}-\mathrm{I})z
\rig +2\lef P^{-1/2}z,
e^{tA}x\rig+\|e^{tA}x\|^2_{P^{-1}_t}\Big]\Big)
\end{aligned}
\end{equation*}
provide
\begin{equation*}
\begin{aligned}
& \int_{\mathbb{R}^d}\big|p(x,t,y)-p(y)\big|dy=
\int_{\mathbb{R}^d}\Big|1-\frac{p(y)}{p(x,t,y)}\Big|p(x,t,y)dy
\\
&=\int_{\mathbb{R}^d}\Bigg|1-\frac{p(P^{1/2}_tz+e^{tA}x)}
{p(x,t,P^{1/2}_tz+e^{tA}x)}\Bigg|p(z)dz
\\
&\le\Bigg|\sqrt{\frac{\det P_t}{\det P}}-1\Bigg| +\sqrt{\frac{\det
P_t}{\det P}}\int_{\mathbb{R}^d}
\Big|\exp\Big(-\frac{1}{2}\Big[\lef
z,(P_tP^{-1}-\mathrm{I})z \rig
\\
&\quad
+2\lef P^{-1/2}z,
e^{tA}x\rig+\|e^{tA}x\|^2_{P^{-1}_t}\Big]\Big)-1\Big|p(z)dz.
\end{aligned}
\end{equation*}
Due to
$\mathbf{(A)}$, $e^{tA}x$ converges to zero in $t\to\infty$
exponentially fast in a sense that $|e^{tA}x|\le ce^{-t\lambda}\|x\|$
for some generic $\lambda>0$.
Moreover, $|P_tP^{-1}-\mathrm{I}|\le \mathfrak{c}e^{-t\lambda}$, owing to $P-P_t$ solves
the differential equation
$
\dot{\triangle}_t=A^*\triangle_t+\triangle_tA
$
subject to $\triangle_0=P$ (see, \eqref{Lyapunov} and \eqref{Lyapt}) .
The above-mentioned convergence implies also
$$
 \Big|\Big(\frac{\det P_t}{\det P}\Big)^{1/2}-1\Big|\le \mathbf{c}e^{-t\lambda}.
$$

Thus, there exists an
appropriate positive continuous function $\upsilon(x) \ (<\infty)$
over $\mathbb{R}^d$ such that for $t\ge t_0>0$,
\begin{gather*}
\int_{\mathbb{R}^d}\big|p(x,t,y)-p(y)\big|dy
\le ce^{-t\lambda}\Big[1+
\int_{\mathbb{R}^d}\{\|z\|^2+\|x\|^2\}e^{\mathfrak{c}e^{-t\lambda}[\|z\|^2+\|x\|^2]}p(z)dz
\Big]
\\
\le ce^{-t\lambda}(1+\upsilon(\|x\|))
\end{gather*}
and, in turn, \eqref{4.3<} holds true, owing to
$$
\int_0^\infty |E_xH(X_s)|ds\le 2Kt_0+\int_{t_0}^\infty|E_xH(X_s)|ds \le
2K+\frac{Kc}{\lambda}c(1+\upsilon(\|x\|)).
$$

\medskip
\underline{Step 3. $\nabla U$ existence.}
Assume $\mathbf{(A''_H)}_{1)}$ and notice that
\begin{equation}\label{nnn}
\int_0^\infty\Big|\int_{\mathbb{R}^d} \nabla_x
H(P_tz+e^{tA}x)\frac{1}{(2\pi)^{d/2}}e^{-\frac{1}{2}\|z\|^2}dz\Big|dt\le \text{const.}
\end{equation}
Since
$
U(x)=-\int_0^\infty\int_{\mathbb{R}^d}
H(P_tz+e^{tA}x)\frac{1}{(2\pi)^{d/2}}e^{-\frac{1}{2}\|z\|^2}dzdt,
$
by virtue of of \eqref{nnn} we have
$$
\nabla U(x)=-\int_0^\infty\int_{\mathbb{R}^d} \nabla_x
H(P_tz+e^{tA}x)\frac{1}{(2\pi)^{d/2}}e^{-\frac{1}{2}\|z\|^2}dzdt.
$$

In particular, $\nabla U$ is bounded.

\smallskip
Assume $\mathbf{(A''_H)}_{2)}$. Now, we prove that
\begin{equation}\label{nnnn}
\int_0^\infty\Big|\int_{\mathbb{R}^d}H(y)\nabla_x p(x,t,y)dy\Big|dt\le \text{const.}
\end{equation}
The use of
$
\nabla_xp(x,t,y)=-p(x,t,y)\big(y-e^{tA}x\big)^*P^{-1}_te^{tA}, \ t>0,
$
provides
\begin{equation*}
\begin{aligned}
\int_{\mathbb{R}^d}H(y)\nabla_x p(x,t,y)dy &=-
EH(X^x_t)(X^x_t-EX^x_t)^*P^{-1}_te^{tA}.
\\
&=E[H(X^x_t)-EH(X^x_t)][X^x_t-EX^x_t]^*P^{-1}_te^{tA}.
\end{aligned}
\end{equation*}
Consequently, taking into account the boundedness of $H$ and \eqref{ebya},
by Cauchy-Schwarz's inequality we get (with a generic positive constant $\lambda$):
\begin{gather*}
\Big|\int_{\mathbb{R}^d}H(y)\nabla_x p(x,t,y)dy\Big|\le
c\Big(E\|X^x_t-EX^x_t\|^2 \Big)^{1/2}|P^{-1}_t||e^tA|
\\
\le \mathfrak{c}e^{-t\lambda}\big(\trace(P_t)\big)^{1/2}|P^{-1}_t|
\le \mathbf{c}\frac{e^{-t\lambda}}{\sqrt{t}}
\end{gather*}

Then, \eqref{nnnn} holds and $\nabla U$ is bounded.

\medskip
\underline{Step 4. $M_t$ existence.} Since an applicability of It\^o's (Krylov-It\^o's)
formula to
$U(X_t)$ is questionable, we show that
$(M_t,\mathscr{F}^X_t)_{t\ge 0}$, with
$$
M_t=U(X_t)-U(x)+\int_0^tH(X_s)ds,
$$
is the continuous martingale,
\begin{equation}\label{new<M>}
\langle M\rangle_t=\int_0^t\nabla^*U(X_s)BB^*\nabla U(X_s)ds
\end{equation}
and $E\|M_t\|^2<\infty$ over $t\in \mathbb{R}_+$; the latter is provided by the
boundedness of $\nabla U$ .

The use of a homogeneity in $t$ of the Markov process $X_t$ enables to claim that
$U(X_t)$ admits the following presentation a.s.,
$$
U(X_t)=\int_t^\infty
E_{X_t}U(X_s)ds=\int_t^\infty E\big(H(X_s)|\mathscr{F}^X_t\big)ds.
$$
Then for any $t'<t$, we have
$$
\begin{aligned}
M_t-M_{t'}&=\int_{t'}^\infty
E\big(H(X_s)|\mathscr{F}^X_t\big)ds-\int_{t'}^\infty
E\big(H(X_s)|\mathscr{F}^X_{t'}\big)ds
\\
& \quad
+\int_{t'}^tE\big(H(X_s)|\mathscr{F}^X_t\big)ds-\int_{t'}^tH(X_s)ds
\quad\text{a.s.}
\end{aligned}
$$
and the martingale property,
$
E(M_t|\mathscr{F}^X_{t'})=M_{t'}
$
a.s., becomes obvious.

Now, we establish \eqref{new<M>} with the help of well known fact:
for any $t>0$, $\langle M\rangle_t$
coincides with the limit, in probability, in $k\to\infty$ of
$$
\sum_{1\le j\le
k}\big(M_{t^k_j}-M_{t^k_{j-1}}\big)\big(M_{t^k_j}-M_{t^k_{j-1}}\big)^*,
$$
where $0\equiv t^k_0<t^k_1<\ldots<t^k_{t_k}\equiv t$ is a
condensing sequence of time values. We recall only that
$
M_{t^k_j}-M_{t^k_{j-1}}=U(X_{t^k_j})-U(X_{t^k_{j-1}})+O(t^k_j-t^k_{j-1})
$
and
$$
U(X_{t^k_j})-U(X_{t^k_{j-1}})=\nabla^*U(X_{t^k_{j-1}})B\big[W_{t^k_j}
-W_{t^k_{j-1}}\big]+O(t^k_j-t^k_{j-1}).
$$

\medskip
\underline{Step 5. (i) verification.} Due to the linear growth condition of
$\|U(x)\|$, suffice it to show that
\begin{equation}\label{4.8i}
\lim_{t\to\infty}\varrho(t)\log P\big(V(X_t)> t^2\varepsilon
\big)=-\infty.
\end{equation}
for $V(x)=\lef x,\Gamma x\rig$ with an appropriate positive definite
matrix $\Gamma$.
In view of $\mathbf{(A)}$, it is convenient to choose $\Gamma$ solving the Lyapunov
equation
$
A^*\Gamma+\Gamma A+\mathrm{I}=0.
$
The function $V(x)$ belongs to the range of definition for
$\mathscr{L}$ with
\begin{equation*}
\begin{aligned}
\mathscr{L}V(x)&=\lef x,(A^*\Gamma+\Gamma A)x\rig +\trace(B\Gamma
B^*)
\\
&=-\|x\|^2+\trace(B\Gamma B^*)\le -cV(x)+\mathfrak{c}
\end{aligned}
\end{equation*}
while
$
V(X_t)-V(x)-\int_0^t\mathscr{L}V(X_s)ds=\int_0^t2\lef X_s,\Gamma BdW_s \rig
=:N_t
$
is the martingale (relative to $(\mathscr{F}^X_t)$) with
$
\langle N\rangle_t=\int_0^t 4\lef X_s,\Gamma^2 X_s\rig dt
\le \int_0^t cV(X_s)ds.
$

Now, \eqref{4.8i} is provided by Corollary \ref{cor-1} to Theorem
\ref{theo-a.1}.

\medskip
\underline{Step 6. (ii) verification.} Since $\nabla U$ is bounded and continuous,
(ii) holds true if
\begin{eqnarray*}
\lim_{t\to\infty}\varrho(t)\log P\Big(\Big|\int_0^th(X_s)ds\Big|>
t\varepsilon\Big)=-\infty \label{4.6d}
\end{eqnarray*}
for any bounded and continuous $h:\mathbb{R}^d\Rightarrow
\mathbb{R}$ with $\int_{\mathbb{R}^d}h(z)p(z)=0$.

Assume for a moment that $h$ satisfy $\mathbf{(A''_H)}_{1)}$ from
Theorem \ref{theo-2}. Then, the function $u(x)=-\int_0^\infty
Eh(X_t)dt$ is well defined and $(u(X_t),\mathscr{F}^X_t)_{t\ge 0}$
is the semimartingale: $ u(X_t)=u(x)-\int_0^th(X_s)ds+m_t, $ where
$(m_t,\mathscr{F}^X_t)_{t\ge 0}$ is the continuous martingale with
$\langle m\rangle_t=\int_0^t\nabla^* u(X_s) BB^*\nabla U(X_s)ds$
and $\nabla u(x)$ is bounded and continuous.

Hence, suffice it to verify
\eqref{4.6d} with $\int_0^th(X_s)ds$ replaced by $u(X_t)$ and $m_t$
separately.

First of all notice that the version of \eqref{4.6d} with $m_t$
is valid due to Theorem \ref{theo-a.2} owing to $\langle m\rangle_t\le Kt$,
where $K\ge \nabla^*u(X_t)BB^*\nabla u(X_t)$ in $t$ over $\mathbb{R}_+$.
Further, because of
$\nabla u$ is bounded and, then, $u$ satisfies the linear growth condition,
the version of \eqref{4.6d} with $u(X_t)$ is reduced to \eqref{4.8i}.

If $h$ does not satisfy $\mathbf{(A''_H)}_{1)}$, we apply the decomposition
$h=h'+h''$ borrowed from the proof of Theorem \ref{theo-V1},
$\alpha=1$. Then, the version of \eqref{4.6d} with $h''$ is reduced to:
for sufficiently large $l$,
\begin{equation*}
\lim_{t\to\infty}\varrho(t)\log P\Big(
\int_0^t\lef X_s,\Gamma X_s\rig ds>t(l^2\varepsilon)\Big)=-\infty,
\end{equation*}
and is verified with the help of Theorem \ref{theo-a.1} for $V(x)=\lef x,\Gamma x\rig$.

The verification of \eqref{4.6d} with $h'$
differs from the corresponding part of proof for Theorem
\ref{theo-V1}, $\alpha=1$. Let $l$, involved in the definition of
$h'$, and $\varepsilon>0$ be chosen. Since $h'$ is compactly
supported, there exists a polynomial $h_\varepsilon$ such that
\begin{equation*}
\begin{aligned}
& c_\varepsilon:=\sup_x|h'(x)-h_\varepsilon(x)|=o(\varepsilon)
\\
&
d_\varepsilon:=\int_{\mathbb{R}^d}h_\varepsilon(z)p(z)dz=o(\varepsilon).
\end{aligned}
\end{equation*}
Because of $\widehat{h}_\varepsilon=h_\varepsilon-d_\varepsilon$ satisfies
\eqref{nul} and $\mathbf{(A''_H)}_{1)}$, the validity of
\eqref{4.6d} with $\widehat{h}_\varepsilon$ is obvious. So, it is left to recall only
that
$\sup_x|h'(x)-\widehat{h}_\varepsilon(x)|=o(\varepsilon)$.
\qed

\subsection{The proof of Theorem \ref{theo-3.4}}
\label{sec-4.4}

Obviously, $H(x)$ satisfies \eqref{nul}.

We shall verify (i), (ii) from Theorem \ref{theo-main}. By virtue
of  \eqref{1.5U}, the quadratic form of $H$ is inherited by $U$.
We examine the following $U(x)=\lef x,\Upsilon x\rig-\upsilon$
with a positive definite matrix $\Upsilon$ and positive number
$\upsilon$.  By It\^o's formula we find that
\begin{equation*}
dU(X_t)=\underbrace{\big[\lef X_t,[\Upsilon A+A^*\Upsilon]X_t\rig
+\trace (B^*\Upsilon B)\big]}_{\text{candidate to be $-H(X_t)$}}dt+
\underbrace{2\lef X_t,\Upsilon BdW_t\rig}_{=M_t}.
\end{equation*}
The realization of this project requires for $\Upsilon$ to be a solution
of Lyapunov's equation
$
\Upsilon A+A^*\Upsilon+\Gamma=0
$
what, in particular, provides
$
\trace (B^*\Upsilon B)=\trace(\Gamma^{1/2}P\Gamma^{1/2}),
$
where $P$ is the covariance of the invariant measure.
With chosen $\Upsilon$, set $\mathsf{D}=\Upsilon BPB^*\Upsilon$ and notice that
$
\langle M\rangle_t=\int_0^t4\lef X_s,\mathsf{D}X_s\rig ds
$

\medskip
(i) is reduced to
$$
\lim_{t\to\infty}\varrho(t)\log P\Big(\lef
X_t,\Gamma X_t\rig>t^\kappa \varepsilon\Big)=-\infty
$$
which holds since for positive and sufficiently small $\lambda$ the moment
generating function
$
\log Ee^{\lambda \lef X_t,\Upsilon X_t\rig }
$
is bounded over $t\in \mathbb{R}_+$ and, then, Chernoff's inequality
provides
$$
\frac{1}{t^{2\kappa-1}}\log P\big(\lef X_t,\Upsilon X_t\rig
>t^\kappa\varepsilon \big) \le -\lambda
t^{1-\kappa}\varepsilon+\frac{\log Ee^{\lambda \lef X_t,\Upsilon
X_t\rig }} {t^{2\kappa-1}} \xrightarrow[t\to\infty]{} -\infty.
$$

\medskip
(ii) is valid if
\begin{equation}\label{4.16}
\lim_{t\to\infty}\varrho(t)\log
P\Big(\Big|\int_0^t\Big[\lef X_s,\Upsilon BB^*\Upsilon
X_s\rig-\trace(\mathsf{D})\Big]ds\Big|
>t\varepsilon\Big)=-\infty.
\end{equation}
Let us denote $\gamma=\Upsilon BB^*\Upsilon$ and
$
\mathfrak{h}(x)=\lef x,\gamma x\rig-\trace(\mathsf{D}).
$
We repeat the previous arguments to find $\mathfrak{u}(x)=\lef
x,\mathfrak{r}x\rig-r$ with a positive definite matrix $\mathfrak{r}$ and
positive number $r$ such that
$
\mathfrak{m}_t=\mathfrak{u}(X_t)-\mathfrak{u}(x)+\int_0^t\mathfrak{h}(X_s)ds
$
is a continuous martingale with $ \langle
\mathfrak{m}\rangle_t=\int_0^t\lef X_s,\mathfrak{q}X_s\rig ds, $
where $\mathfrak{q}$ is a positive definite matrix. Now, we may replace
\eqref{4.16} by
\begin{equation*}
\begin{aligned}
& {\rm (1)} \ \lim_{t\to\infty}\varrho(t)\log P\big(\lef
X_t,\gamma X_t\rig
>t\varepsilon\big)=-\infty
\\
& {\rm (2)} \ \lim_{t\to\infty}\varrho(t)\log
P\big(\big|\mathfrak{m}_t\big|
>t\varepsilon\big)=-\infty.
\end{aligned}
\end{equation*}

(1) is verified similarly to (i).
(2) is verified with the help of Theorem \ref{theo-a.2} by showing
$
\lim_{t\to\infty}\varrho(t)\log P\big(\langle \mathfrak{m}\rangle_t> tn\big)=-\infty
$
for sufficiently large $n$ what is nothing but
\begin{equation}\label{3.10aa}
\lim_{t\to\infty}\varrho(t)\log P\Big(\int_0^t\lef X_s,\mathfrak{q}X_s\rig
ds> tn\Big)=-\infty.
\end{equation}
A version of \eqref{3.10aa} with $\mathfrak{q}$ replaced by any positive definite
matrix $G$ provides  \eqref{3.10aa} too. For computational convenience,
we take  $G$ solving Lyapunov's equation
$ A^*G+GA+\mathrm{I}=0.$ The function
$V(x)=\lef x,Gx\rig$ belongs to the range of definition of $\mathscr{L}$
with
$
\mathscr{L}V(x)=-2\|x\|^2+\trace(BPB^*)\le -cV(x)+\mathfrak{c}
$
\ and
$$
N_t=V(X_t)-V(x)-\int_0^t\mathscr{L}V(X_s)ds=
\int_0^t2\lef X_s,BdW_s\rig,
$$
so that,
$
\langle N\rangle_t\le\int_0^t\mathbf{c}V(X_s)ds.
$

Thus, the proof is completed by applying
Theorem \ref{theo-a.1}. \qed

\section{Example of statistical application}
\label{sec-5}

Let $X_t \ (\in\mathbb{R})$ be a diffusion process: $ dX_t=-\theta
X_tdt+dW_t, $ subject to a fixed $X_0$. The parameter $\theta\in
(0,\infty)$ is unknown and is evaluated with help of well known estimate
\begin{equation*}
\widehat{\theta}_t=\frac{\int_0^tX_sdX_s}{\int_0^tX^2_sds}, \ t>0.
\end{equation*}
It is well known that the CLT holds for the family
$\big(\sqrt{t}(\theta-\widehat{\theta}_t)\big)_{t\to\infty}$ with a limit: zero mean
Gaussian random variable with the  variance $2\theta$.

In this section, we show that $\theta-\widehat{\theta}_t$ possesses an asymptotic
(in $t\to\infty$) in the MDP scale, $\frac{1}{2}<\kappa<1$, that is,
the family
$\big(t^{1-\kappa}(\theta-\widehat{\theta}_t)\big)_{t\to\infty}$
obeys $(\varrho,J)$-MDP with $J(Y)=\frac{Y^2}{4\theta}$.
The use of some details from the proof of Theorem \ref{theo-3.4} enables to claim that
$\big(\frac{1}{t}\int_0^t[X^2_s-\frac{1}{2\theta}]ds)_{t\to\infty}$ is negligible
in $\varrho$-MDP scale. Therefore, the family
$\big(t^{1-\kappa}(\theta-\widehat{\theta}_t)\big)_{t\to\infty}$
shares the MDP with
$
\big(\frac{1}{t^\kappa}\int_0^t2\theta X_sdW_s\big)_{t\to\infty}.
$
Further, the announced MDP hold if (ii) from Theorem \ref{theo-main} is valid:
\begin{equation}\label{5.2}
\lim_{t\to\infty}\log\varrho(t) P\Big(\Big|\int_0^t[4\theta^2X^2_s-Q]ds\Big|
>t\varepsilon\Big)=-\infty.
\end{equation}
Obviously, $Q=2\theta$ and the validity of \eqref{5.2} is verified with the help of
arguments used in the proof of Theorem \ref{theo-3.4}.

In particular, this MDP and the contraction Varadhan's principle, for sufficiently
large $t$ provide
\begin{equation*}
\frac{1}{t^{2\kappa-1}}\log P\big(t^{1-\kappa}
|\theta-\widehat{\theta}_t|>\delta\big)\asymp -\frac{\delta^2}{4\theta}.
\end{equation*}

\appendix
\section{Exponential negligibility of functionals and martingales}
\label{sec-A}

Let $X_t$ be a diffusion process defined in \eqref{ito00} with $X_0=x$.

Assume $V(x):\mathbb{R}^d\to\mathbb{R_+}$, with
$\lim_{\|x\|\to\infty}V(x)=\infty$, belongs to
the range of definition of $\mathscr{L}$. Introduce a martingale
relative to $(\mathscr{F})_{t\ge 0}$:
\begin{equation}\label{aa}
N_t=V(X_t)-V(x)-\int_0^t\mathscr{L}V(X_s)ds.
\end{equation}

\begin{theorem}\label{theo-a.1}
Assume

{\rm 1)} $\mathscr{L}V\le -cV^\ell+\mathfrak{c}$, $\exists \ \ell>0$

{\rm 2)} $\langle N\rangle_t\le
\int_0^t\mathbf{c}\big(1+V^r(X_s)\big)ds$, $\exists \ r\le \ell$.

\noindent Then, for any $\varepsilon>0$ and sufficiently large
number $n$
\begin{eqnarray*}
&& \lim\limits_{t\to\infty}\varrho(t)\log
P\big(V(X_t)>t^{2\kappa}\varepsilon\big) =-\infty,
\\
&& \lim\limits_{t\to\infty}\varrho(t)\log P\Big(\int_0^tV^\ell(X_s)ds
>tn\Big)
=-\infty.
\end{eqnarray*}
over $x\in\mathbb{R}^d$.
\end{theorem}
\begin{corollary}
\label{cor-1}
$
\lim\limits_{t\to\infty}\varrho(t)\log
P\big(V(X_t)>t^2\varepsilon\big) =-\infty,
$
since $t^2>t^{2\kappa}, t>1.$
\end{corollary}
\begin{remark}\label{rem-4}
The statements of Theorem \ref{theo-a.1} remain valid if constants
$c, \mathfrak{c}, \mathbf{c}$, involved in {\rm 1)} and {\rm 2)}
depend on $\varepsilon$.
\end{remark}

\begin{theorem}\label{theo-a.2}
Let $M_t(\in\mathbb{R}, M_0=0)$ be a continuous martingale.

Then, for any $\varepsilon>0$,
$$
\lim_{t\to\infty}\varrho(t)\log
P\big(|M_t|>t\varepsilon\big)=-\infty
$$
provided that, under sufficiently large number $n$ depending on $\varepsilon$,
$$
\lim_{t\to\infty}\varrho(t)\log P\big(\langle
M\rangle_t>tn\big)=-\infty.
$$
\end{theorem}

\medskip
\noindent {\it The proof of Theorem \ref{theo-a.1}.} With
$\lambda\in\mathbb{R}$, and the (continuous) martingale $N_t$ from \eqref{aa},
we introduce a positive random process
$
\mathfrak{z}_t(\lambda)=e^{\lambda N_t-0.5\lambda^2\langle
N\rangle_t}.
$
It is well known and easily verified with the help of It\^o's formula that
$(\mathfrak{z}_t(\lambda),\mathscr{F}^X_t)_{t\ge 0}$ is
a positive local martingale. Moreover, by Problem 1.4.4, \cite{LSMar},
it is a supermartingale too. We shall use the supermartingale property:
 $E\mathfrak{z}_t(\lambda)\le E\mathfrak{z}_0(\lambda)\equiv 1$
over $t\in\mathbb{R}_+$.
Denote by
$$
\mathfrak{A}_1=\{V(X_t)>t^{2\kappa}\varepsilon\}\quad\text{and}\quad
\mathfrak{A}_2=\Big\{\int_0^tV^r(X_s)ds>tn\Big\}.
$$
The use of $E\mathfrak{z}_t(\lambda)\le 1$ provides
\begin{equation}\label{a.6}
1\ge EI_\mathfrak{{A}_i}\mathfrak{z}_t(\lambda), \quad i=1,2
\end{equation}
Notice that \eqref{a.6} remains valid with $\mathfrak{z}_t(\lambda)$ replaced by its
lower bound on $\mathfrak{A}_i$. We proceed the proof by finding appropriate
deterministic (!) lower bounds.
Write
\begin{equation*}
\lambda N_t-0.5\lambda^2\langle
N\rangle_t=\lambda\Big(V(X_t)-V(x)-
\int_0^t\mathscr{L}V(X_s)ds\Big)-0.5\lambda^2\langle N\rangle_t.
\end{equation*}
Thence, in view of 1) and 2), with $\lambda>0$ we get
\begin{equation*}
\begin{aligned}
\lambda N_t-0.5\lambda^2\langle N\rangle_t&\ge
\lambda\Big(V(X_t)-V(x)+ \int_0^t[cV^\ell (X_s)-\mathfrak{c}]ds\Big)
\\
&\quad -0.5\lambda^2\int_0^t\mathbf{c}(1+V^r(X_s))ds.
\end{aligned}
\end{equation*}
Taking into account $1+V^r(X_s)\le 2+V^\ell (X_s)$, provided by $r\le \ell$,
and choosing
$\lambda^\circ=\argmax_{\lambda>0}\big[c\lambda-0.5\mathbf{c}\lambda^2\big]
=\frac{c}{\mathbf{c}}$, we get
\begin{equation*}
\begin{aligned}
\lambda^\circ N_t-0.5(\lambda^\circ)^2\langle N\rangle_t&\ge
\frac{c}{\mathbf{c}}\big[V(X_t)-V(x)\big]-t\frac{c}{\mathbf{c}}\big[\mathfrak{c}+c]
+ \frac{c^2}{2\mathbf{c}}\int_0^tV^\ell(X_s)ds
\\
&\ge
  \begin{cases}
    \frac{c}{\mathbf{c}}\big[t^{2\kappa}\varepsilon-V(x)\big]
    -t\frac{c}{\mathbf{c}}\big[\mathfrak{c}+c], & \text{over $\mathfrak{A}_1$}, \\
    -\frac{c}{\mathbf{c}}V(x)-t\frac{c}{\mathbf{c}}\big[\mathfrak{c}+c]
+ \frac{c^2}{2\mathbf{c}}tn, & \text{over $\mathfrak{A}_2$}.
  \end{cases}
\end{aligned}
\end{equation*}
These lower bounds jointly with \eqref{a.6} provide
$$
\left. \begin{array}{ll}
\varrho(t)\log P\big(\mathfrak{A}_1\big)\le -\frac{c}{\mathbf{c}}
\Big[t\varepsilon-\frac{V(x)}{t^{2\kappa-1}}\Big]
    +t^{2(1-\kappa)}\frac{c}{\mathbf{c}}\big[\mathfrak{c}+c], & \text{over $\mathfrak{A_1}
    $}
    \\
\varrho(t)\log P\big(\mathfrak{A}_2\big)\le
 \frac{c}{\mathbf{c}}\frac{V(x)}{t^{2\kappa-1}}+t^{2(1-\kappa)}\frac{c}{\mathbf{c}}
 \big[\mathfrak{c}+c]
- t^{2(1-\kappa)}\frac{c^2}{2\mathbf{c}}n, & \text{over $\mathfrak{A}_2$}
\end{array}
\right\}\xrightarrow[t\to\infty]{}-\infty.
$$
\qed

\bigskip
\noindent {\it The proof of Theorem \ref{theo-a.2}.}
Notice that only
\begin{equation}\label{aa.3}
\lim_{t\to\infty}\varrho(t)\log P\big(|M_t|>t\varepsilon, \langle
M\rangle_t\le tn\big)=-\infty
\end{equation}
is required to be proved. Moreover, it suffices to prove only
\begin{equation}\label{aa.4}
\lim_{t\to\infty}\varrho(t)\log P\big(M_t>t\varepsilon, \langle
M\rangle_t\le tn\big)=-\infty
\end{equation}
owing to a version with $-M_t$ is verified similarly and both ``$\pm M_t$''
provide \eqref{aa.3}.

For \eqref{aa.4} verification, we use the inequality from
\eqref{a.6} with $\lambda>0$ and
$N_t$, $\langle N\rangle_t$ replaced by $M_t$, $\langle M\rangle_t$
respectively and $\mathfrak{A}_i$ replaced by
$$
\mathfrak{A}=\mathfrak\{M_t>t\varepsilon,\langle M\rangle_t\le
tn\}
$$
and notice that
$$
\log \mathfrak{z}_t(\lambda)=\lambda M_t-0.5\lambda^2\langle
M\rangle_t\underbrace{\ge}_{\text{over }\mathfrak{A}} \lambda
t\varepsilon-0.5\lambda^2 tn\ge \min_{\lambda>0} (\lambda
t\varepsilon-0.5\lambda^2 tn)=t\frac{\varepsilon^2}{2n}.
$$
Then, owing to $1\ge e^{t\frac{\varepsilon^2}{2n}}EI_{\mathfrak{A}}$, we get
$ \varrho(t)\log P\big(\mathfrak{A}\big)\le
-t^{2(1-\kappa)} \frac{\varepsilon^2}{2n}\to-\infty. $ \qed

\end{document}